\documentclass{l4dc2027}

\title[Finite sample guarantees of MPPI]{Finite sample guarantees of MPPI: suboptimality, sample complexity and closed-loop guarantees}
\usepackage{times}

\usepackage{amssymb}  
\usepackage{mathrsfs}
\usepackage{mathtools}
\usepackage{caption}

\usepackage{enumitem}
\setlist[enumerate]{left=0pt, label=(\arabic*)}
\setlist[itemize]{left=0pt}
\usepackage{xcolor}

\newtheorem{assumption}{Assumption}

\coltauthor{%
  \Name{Changyi Lei} \Email{chalei1@student.unimelb.edu.au}\\
  \Name{Dragan Ne\v{s}i\'{c}} \Email{dnesic@unimelb.edu.au}\\
  \Name{Ye Pu} \Email{ye.pu@unimelb.edu.au}\\
  \addr Department of Electrical and Electronic Engineering,
  University of Melbourne, Parkville, VIC 3010, Australia
}

\begin{document}

\maketitle

\begin{abstract}
This paper considers model predictive path integral (MPPI) control for nonlinear finite-horizon optimal control problems with bounded inputs, allowing nonconvex and nonsmooth costs. We prove theoretical guarantees that quantify the suboptimality caused by positive temperature and finite sampling on MPPI and establish their implications for closed-loop guarantees. Specifically, we derive a suboptimality bound for infinite-sample MPPI that characterizes the effects of temperature and prior coverage, followed by a high-probability finite-sample approximation error bound that reveals the dependence on sample size, temperature, and cost variability. By incorporating these bounds into a terminal-constraint-free MPC framework, we further establish semi-global practical recurrence of the closed-loop system. Numerical examples involving a four-link planar manipulator illustrate the temperature–sampling tradeoff and the influence of prior quality on approximation accuracy and closed-loop performance.
\end{abstract}


\section{Introduction}

Model predictive control (MPC) is widely used across engineering applications due to its ability to optimize predicted performance while explicitly accounting for constraints \citep{borrelli2017predictive}. However, applying MPC to complex dynamical systems remains challenging, particularly when nonlinear dynamics, nonconvex objectives, and nonsmooth phenomena such as contact and hybrid transitions arise. Conventional nonlinear optimization methods, including SQP \citep{1470308}, interior-point methods \citep{wachter2006implementation}, and iLQR/DDP \citep{li2004iterative,tassa2012synthesis}, typically exploit derivatives or local smooth approximations. While effective for many problems, their performance can be sensitive to nonconvexity and nonsmoothness, motivating optimization approaches with fewer structural requirements.

Among these approaches, model predictive path integral (MPPI) control has attracted increasing attention \citep{honda2026model}. MPPI evaluates sampled control sequences through forward simulations, assigns importance weights based on their trajectory costs, and computes a cost-weighted control update \citep{williams2017model}. By relying on trajectory evaluations rather than derivatives, MPPI offers flexibility for nonlinear, nonconvex, and nonsmooth optimal control problems without requiring explicit differentiation through the dynamics. Moreover, its independent rollouts are naturally parallelizable, enabling efficient implementation on GPUs \citep{williams2017model}. These advantages have motivated applications in autonomous driving \citep{williams2018information}, aerial vehicles \citep{minavrik2024model, 11386891}, trajectory planning \citep{chen2026cablerobotgraphsim,pmlr-v331-hori26a} and locomotion \citep{xue2025full, 11696779}, as well as numerous methodological extensions \citep{honda2024stein, 9993190, 11312129, 10161100, 10824881, 11373877}.



Despite its empirical success, existing theoretical guarantees for MPPI remain fragmented. Existing analyses can be broadly divided into two directions with distinct limitations. (1) The first studies MPPI from an optimization perspective assuming infinite samples. CoVO-MPC \citep{yi2024covo} establishes linear contraction for infinite-sample MPPI under smooth and strongly convex costs. \citep{homburger2025optimality} analyzes infinite-sample MPPI in a smooth, unconstrained setting with a unique optimizer, showing asymptotic first- and second-order control and value suboptimality in the temperature as it vanishes. \citep{fazlyab2026model} interprets MPPI as preconditioned gradient descent on a
KL-regularized cost and establishes descent and stationarity
for the infinite-sample MPPI. \citep{wang2026generalized} interprets MPPI as an expectation--maximization (EM) algorithm and establishes asymptotic convergence of repeated, infinite-sample MPPI updates to the optimum. While these studies provide important insights into the optimization mechanisms of MPPI, they generally rely on smoothness or other structural conditions and do not provide finite-sample error bound. The second direction studies finite-sample approximation error and closed-loop stability. \citep{yoon2022sampling} derives probabilistic sample-complexity bounds that quantify the finite-sample Monte Carlo approximation error of infinite-sample MPPI. Similar analysis is used in \citep{yoon2026finite} to help establish finite-sample practical stability by viewing MPPI as a stochastic perturbation of the stabilizing LQR controller. The stability analysis is further extended to nonlinear systems in \citep{yoon2026stochastic}. However, existing stability results either rely on specialized problem structures or involve nonconstructive constants that are difficult to quantify, limiting the derivation of explicit stability conditions for general nonlinear systems. Consequently, a unified nonasymptotic analysis that quantifies both temperature-induced suboptimality and finite-sample approximation error, and connects these bounds to closed-loop guarantees under general nonconvex and nonsmooth costs, remains lacking.

Motivated by these gaps, this paper develops a unified positive-temperature and finite-sample analysis of MPPI that directly connects suboptimality, sample complexity, and closed-loop guarantees beyond the restrictive convex or smooth settings considered in existing works. Our contributions are:
\begin{enumerate}
    \item An explicit positive-temperature suboptimality bound is established for infinite-sample MPPI relative to the original finite-horizon optimum under only local Lipschitz continuity, without requiring convexity or differentiability. General sampling distributions and convex bounded input constraints are accommodated, with the effects of temperature and prior distribution quantified explicitly in a form suitable for closed-loop stability analysis.

    \item A high-probability finite-sample approximation bound is derived, with a transparent characterization of the separate effects of temperature, sample size, and cost variability.

    \item The positive-temperature and finite-sample errors are propagated through a terminal-constraint-free MPC framework to obtain a semi-global practical recurrence guarantees in the sense of \citep{GRAMMATICO20132939,moldenhauer2025robust}.
\end{enumerate}

\textit{Notation}: 
For two probability measures $Q$ and $P$, the notation $Q\ll P$ means that $Q$ is absolutely continuous with respect to $P$, i.e., $P(A)=0$ implies $Q(A)=0$ for every measurable set $A$. For two probability measures \(Q\ll P\), we denote the chi-squared divergence by
\(
    \chi^2(Q\|P)
    :=
    \int
    \left(
        \frac{dQ}{dP}-1
    \right)^2 dP.
\)
We use
\(
    \Pi_0
    :=
    \begin{bmatrix}
        I_m & 0_{m\times m(N-1)}
    \end{bmatrix}
    \in\mathbb R^{m\times mN}
\)
to extract the first control input from the stacked sequence
\(\mathbf u=\operatorname{col}(u_0,\ldots,u_{N-1})\), i.e.,
\(\Pi_0\mathbf u=u_0\).
Let
\(
    B_a(\mathbf u^\star)
    :=
    \left\{
        \mathbf u\in\mathbb R^{mN}:
        \|\mathbf u-\mathbf u^\star\|\le a
    \right\}
\)
denote the Euclidean ball of radius \(a>0\) centered at
\(\mathbf u^\star\). For a set \(\mathcal A\subseteq\mathbb R^d\),
\(
    \operatorname{TN}_{\mathcal A}(m,\Sigma)
\)
denotes the Gaussian distribution \(\mathcal N(m,\Sigma)\) truncated to
\(\mathcal A\) with density being
\(
    p(z)
    =
    \frac{
        \varphi_d(z;m,\Sigma)\mathbf 1_{\mathcal A}(z)
    }{
        \int_{\mathcal A}\varphi_d(v;m,\Sigma)\,dv
    },
\)
where \(\varphi_d(\cdot;m,\Sigma)\) is the density of \(\mathcal N(m,\Sigma)\). Let $\mathbf 1_{\mathcal A}(z)$ denotes the indicator function, which equals $1$ if $z\in\mathcal A$ and $0$ otherwise. For a square-integrable random variable $X$, we denote its variance by
$\operatorname{Var}(X):=\mathbb E[(X-\mathbb E[X])^2]$. We denote by $\mathcal{K}$ the class of continuous, strictly increasing
functions $\alpha:\mathbb{R}_{\geq 0}\to\mathbb{R}_{\geq 0}$ satisfying
$\alpha(0)=0$, by $\mathcal{K}_\infty$ the subclass of $\mathcal{K}$-functions satisfying $\lim_{s\to\infty}\alpha(s)=\infty$, and by $\mathcal{KL}$ the class of functions
$\beta:\mathbb{R}_{\geq 0}\times\mathbb{R}_{\geq 0}\to\mathbb{R}_{\geq 0}$
such that $\beta(\cdot,t)\in\mathcal{K}$ for each fixed $t\geq 0$ and
$\beta(s,t)$ decreases to zero as $t\to\infty$ for each fixed $s>0$.


\section{Problem Setting and Algorithm}
\label{sec:problem_setting}

\subsection{System Dynamics and Optimization Problem}

Consider the deterministic nonlinear discrete-time system
\begin{equation}
    x_{t+1}=f(x_t,u_t),
    \, t\in\mathbb N_0,
    \label{eq:system}
\end{equation}
where
\(x_t\in \mathbb R^n\)
denotes the state,
\(u_t\in\mathbb U\subset\mathbb R^m\)
denotes the control input and $f:\mathbb R^n \times \mathbb U \rightarrow \mathbb R^n$ is the dynamics that may be nonsmooth. The input constraint set $\mathbb U\subset\mathbb R^m$ is assumed to be nonempty, compact, and convex. Moreover, there exists a constant $u_{\max}>0$ such that
\(
    \mathbb U
    \subseteq
    \left\{
        u\in\mathbb R^m:
        \|u\|_\infty\leq u_{\max}
    \right\}.
\)
For a prediction horizon \(N\in\mathbb N\), define the stacked control sequence
\(
    \mathbf u
    :=
    [u_0^\top,\ldots,u_{N-1}^\top]^\top
    \in
    \mathbb U^N
    \subset
    \mathbb R^{mN}.
\)


Let
\(
\ell:\mathbb R^n \times\mathbb U\rightarrow\mathbb R_{\geq0}
\)
and
\(
\ell_f:\mathbb R^n\rightarrow\mathbb R_{\geq0}
\)
denote the stage cost and terminal cost, respectively. For a given state \(x\), the performance of a candidate control sequence \(\mathbf u\) over the prediction horizon is measured by the finite-horizon cost
\begin{equation}
\begin{split}
    J_N(x,\mathbf u)
    &:=
    \sum\nolimits_{k=0}^{N-1}
    \ell\!\left(
        \hat{x}_k,u_k
    \right)
    +
    \ell_f\!\left(
        \hat{x}_N
    \right), \\
    & s.t. \ \hat{x}_0=x, \ \hat{x}_{k+1}=f(\hat{x}_k, u_k), \ u_k \in \mathbb U, \ k=0,1,...,N-1.
    \end{split}
    \label{eq:finite_horizon_cost}
\end{equation}
The control objective is to select a sequence of admissible inputs that
minimizes \eqref{eq:finite_horizon_cost} while satisfying the
input constraints along the predicted trajectory. For any \(x\in\mathbb R^n\), the optimal value function and the optimal control sequence are respectively
\begin{equation}
    V_N(x)
    :=
    \min_{\mathbf u\in \mathbb U^N}
    J_N(x,\mathbf u), \ \mathbf u^\star(x)
    :=
    \arg\min_{\mathbf u\in \mathbb U^N}
    J_N(x,\mathbf u).
    \label{eq:optimal_value_function}
\end{equation}

The following assumptions are imposed.
\begin{assumption}[Bounded finite-horizon cost]
\label{ass:convex_bounded_cost}
For any $N\in\mathbb N$ and $x\in\mathbb R^n$, there exists a finite constant \(C_N(x)>0\) such that
\(
    0
    \leq
    J_N(x,\mathbf u)
    \leq
    C_N(x),
    \,
    \forall
    \mathbf u\in
    \mathbb U^N.
\)
\end{assumption}

\begin{assumption}[Local Lipschitz continuity \citep{liu2026certainty,gandhi2021robust}]
\label{ass:local_lipschitz_cost}
For every $c>0$, the value-function sublevel set
$\Omega_c:=\{x\in\mathbb R^n:V_N(x)\le c\}$
is compact. There exist finite positive constants
$L_{f,c}^x$, $L_{f,c}^u$, $L_{c,N}^x$, and $L_{c,N}^u$
such that, for all $x,y\in\Omega_c$, 
$u,v\in\mathbb U$ and $\mathbf u,\mathbf v\in\mathbb U^N$,
\begin{equation}
    \|f(x,u)-f(y,v)\|
    \le
    L_{f,c}^x\|x-y\|
    +
    L_{f,c}^u\|u-v\|.
    \label{eq:f_local_lipschitz}
\end{equation}
\begin{equation}
    |J_N(x,\mathbf u)-J_N(y,\mathbf v)|
    \le
    L_{c,N}^x\|x-y\|
    +
    L_{c,N}^u\|\mathbf u-\mathbf v\|.
    \label{eq:J_local_lipschitz}
\end{equation}
\end{assumption}

\begin{assumption}[Semi-global finite-horizon cost domination, \citep{grune2012nmpc, grimm2005model}]
\label{ass:no_terminal_cost_condition}
There exists a continuous, positive definite function $\sigma:\mathbb R^n\to\mathbb R_{\geq 0}$ satisfying $\sigma(x)\to\infty$ as $\|x\|\to\infty$ such that
\(
    \ell(x,u)
    \geq
    \sigma(x),
    \,
    \forall x\in\mathbb R^n,
    \,
    \forall u\in\mathbb U.
\)
Moreover, there exists a lcoally bounded function
\(
\bar a:\mathbb R^n\rightarrow\mathbb R_{\geq1}
\)
such that
\(
    V_N(x)
    \leq
    \bar a(x)\sigma(x),
    \,
    \forall x\in\mathbb R^n, \forall N\in\mathbb N.
\)
\end{assumption}

\begin{remark}[Discussion of Assumptions 1--3]
Assumption~1 ensures boundedness of the finite-horizon cost over admissible control sequences, which is used to establish the explicit finite-sample approximation bounds in Section \ref{subsec:sample_complexity_mppi}. This condition is satisfied when the dynamics and cost functions are continuous with bounded inputs. Assumption~2 provides the local Lipschitz properties needed to translate control approximation errors into cost suboptimality and closed-loop stability guarantees. These properties hold when the dynamics, stage cost, and terminal cost are locally Lipschitz continuous on their respective domains. Assumption~3 imposes a finite-horizon cost domination condition, which is used to establish the closed-loop stability guarantees. This condition is commonly employed in MPC stability analysis without terminal constraints.
\end{remark}


\subsection{Finite-Sample MPPI Algorithm}
\label{subsec:finite_sample_mppi}

We consider a standard MPPI controller for approximately solving the finite-horizon optimal control problem \eqref{eq:optimal_value_function}. Rather than directly solving the nonlinear optimization problem \eqref{eq:optimal_value_function}, MPPI draws candidate control sequences from a sampling distribution and assigns larger weights to sequences with lower realized finite-horizon cost.

At any given timestep, let
\(
    \mathbf u^{(1)},\ldots,\mathbf u^{(M)}
    \overset{\mathrm{i.i.d.}}{\sim}P
\)
denote $M$ independent control-sequence samples $\mathbf u^{(m)}\in\mathbb U^N\subset\mathbb R^{mN}$ from the user-defined prior distribution $P\in\mathcal P(\mathbb U^N)$. For each sampled sequence,
define
\(
    J_N^{(m)}(x)
    :=
    J_N\!\left(x,\mathbf u^{(m)}\right),
    \,
    m=1,\ldots,M.
\)
Given a temperature parameter \(\lambda>0\), the finite-sample MPPI control law associates with each
sample the normalized exponential weight
\begin{equation}
    \omega_M^{(m)}(x)
    :=
    \frac{
        \exp\!\left(
            -J_N^{(m)}(x)/\lambda
        \right)
    }{
        \sum\nolimits_{j=1}^{M}
        \exp\!\left(
            -J_N^{(j)}(x)/\lambda
        \right)
    },
    \,
    m=1,\ldots,M.
    \label{eq:mppi_weights}
\end{equation}
The resulting finite-sample MPPI control sequence is the weighted
average
\begin{equation}
    \mathbf u_{M,\lambda}(x)
    :=
    \sum\nolimits_{m=1}^{M}
    \omega_M^{(m)}(x)\mathbf u^{(m)}.
    \label{eq:finite_sample_mppi}
\end{equation}
The receding-horizon implementation is summarized in Algorithm~\ref{alg:mppi}.

\begin{algorithm2e}[t]
\caption{Finite-Sample MPPI Algorithm}
\label{alg:mppi}
\LinesNumbered
\DontPrintSemicolon

Choose horizon $N$, sample size $M$, temperature $\lambda>0$,
and prior distribution $P$\;

\For{$t=0,1,2,\ldots$}{
    Observe $x_t$\;

    \For{$m=1,\ldots,M$}{
        Sample $\mathbf u^{(m)}\sim P$\;

        Evaluate finite-horizon cost
        $J_N^{(m)}(x_t)$\;
    }

    Compute \eqref{eq:mppi_weights} and
    \eqref{eq:finite_sample_mppi}, and apply
    $u_t\gets\Pi_0\mathbf u_{M,\lambda}(x_t)$\;
}
\end{algorithm2e}

\begin{remark}[Input admissibility]
We require the sampling prior $P$ to be supported on
$\mathbb U^N$ as $P\in\mathcal P(\mathbb U^N)$ so that all sampled control sequences are admissible. Since $\mathbb U$ is assumed to be convex, the finite-sample MPPI control sequence $\mathbf u_{M,\lambda}(x)$ remains in $\mathbb U^N$ for any sample size $M$. 
\end{remark}

\section{Main Results}
This section develops the main theoretical results for finite-sample MPPI in Algorithm \ref{alg:mppi}. We first characterize the infinite-sample solution $\mathbf u_{\infty,\lambda}$ through a KL-regularized optimal control problem in Section~\ref{subsec:ideal_mppi}, and then quantify its suboptimality relative to the optimal solution $\mathbf u^\star$ defined as \eqref{eq:optimal_value_function} in Section~\ref{subsec:ideal_mppi_suboptimality}. Building on this result, Section~\ref{subsec:sample_complexity_mppi} establishes a bound on the finite-sample approximation error between the infinite-sample solution $u_{\infty, \lambda}$ and the finite-sample solution $u_{M,\lambda}$ given by Algorithm \ref{alg:mppi}. Finally, combining these two error bounds and building upon MPC without terminal constraints framework, Section~\ref{subsec:stability} establishes semi-global practical recurrence through closed-loop Lyapunov analysis. Omitted proofs are provided in the appendices.

\subsection{Variational Characterization of Infinite-sample MPPI}
\label{subsec:ideal_mppi}

We first characterize the infinite-sample MPPI, which serves as the reference for quantifying the approximation error introduced by finite sampling. For any fixed state $x$, applying the strong law of
large numbers to \eqref{eq:finite_sample_mppi} yields
\begin{align}
\mathbf u_{\infty,\lambda}(x)
&:=\lim_{M\to\infty}\mathbf u_{M,\lambda}(x)
=
\frac{
\mathbb E_{\mathbf u\sim P}
\left[
e^{-J_N(x,\mathbf u)/\lambda}\mathbf u
\right]
}{
\mathbb E_{\mathbf u\sim P}
\left[
e^{-J_N(x,\mathbf u)/\lambda}
\right]
},
\label{eq:ideal_mppi_control}
\end{align}
where the convergence holds almost surely. The following lemma establishes a variational interpretation of this limit through a KL-regularized optimal control problem. This characterization will be instrumental in the subsequent suboptimality analysis.

\begin{lemma}[Variational characterization of infinite-sample MPPI]
\label{lem:ideal_mppi_variational}
For any given $x$, $\lambda>0$ and prior distribution $P$, the infinite-sample
MPPI solution \eqref{eq:ideal_mppi_control} satisfies
\(
\mathbf u_{\infty,\lambda}(x)
=
\mathbb E_{\mathbf u\sim Q_{\lambda,x}^\star}[\mathbf u],
\)
where $Q_{\lambda,x}^\star$ is the optimal gibbs posterior satisfying
\begin{equation}
Q_{\lambda,x}^\star
=
\arg\min_{Q\ll P}
\left\{
\mathbb E_{\mathbf u\sim Q}[J_N(x,\mathbf u)]
+
\lambda\,\mathrm{KL}(Q\|P)
\right\}, \ 
\frac{dQ_{\lambda,x}^\star}{dP}(\mathbf u)
=
\frac{
e^{-J_N(x,\mathbf u)/\lambda}
}{
\mathbb E_{\mathbf u\sim P}
\left[e^{-J_N(x,\mathbf u)/\lambda}\right]
}.
\label{eq:ideal_mppi_posterior}
\end{equation}
\end{lemma}

Thus, infinite-sample MPPI computes the mean of the optimal Gibbs posterior, while finite-sample MPPI approximates this mean through self-normalized Monte Carlo sampling. Importantly, although $Q_{\lambda,x}^\star$ minimizes the KL-regularized objective, its mean $\mathbf u_{\infty,\lambda}$ need not minimize the original finite-horizon cost \eqref{eq:finite_horizon_cost}. This distinction motivates the suboptimality analysis in
Section~\ref{subsec:ideal_mppi_suboptimality}, where we quantify the intrinsic suboptimality of $\mathbf u_{\infty,\lambda}$ relative to $\mathbf u^\star$ defined in \eqref{eq:optimal_value_function}.

We define the infinite- and finite-sample receding-horizon feedback laws, respectively, as
\begin{equation}
\!\kappa_{\infty,\lambda}(x)
:=
\Pi_0\mathbf u_{\infty,\lambda}(x), \,
\kappa_{M,\lambda}(x)
:=
\Pi_0\mathbf u_{M,\lambda}(x), \forall x \in \mathbb R^n.
\label{eq:finite_sample_mppi_feedback}
\end{equation}


\subsection{Suboptimality of Infinite-sample MPPI} \label{subsec:ideal_mppi_suboptimality}
This subsection quantifies the intrinsic suboptimality of the infinite-sample MPPI control law $\mathbf u_{\infty,\lambda}(x)$ defined in \eqref{eq:ideal_mppi_control}, relative to the globally optimal control sequence $\mathbf u^\star(x)$ in \eqref{eq:optimal_value_function}. The main idea is that, when the global minimizer is separated by a positive cost gap from control sequences outside any fixed neighborhood and the prior assigns nonvanishing probability mass to its neighborhood, the posterior distribution $Q_{\lambda,x}^*$ in \eqref{eq:ideal_mppi_posterior} increasingly concentrates around \(\mathbf u^\star(x)\) as the temperature \(\lambda\) decreases. We formalize these conditions to bound \( \|\mathbf u_{\infty,\lambda}(x)-\mathbf u^\star(x) \| \) and the resulting cost suboptimality \(J_N(x,\mathbf u_{\infty,\lambda}(x))-V_N(x)\), borrowing ideas from a nonasymptotic Laplace principle argument \citep[Proposition 4.5]{fornasier2024consensus}.

For any $c>0$, $R>0$ and the value-function sublevel set
$\Omega_c:=\{x\in\mathbb R^n:V_N(x)\le c\}$, define the uniform separation gap
\begin{equation}
\gamma_{N,c}(R)
:=
\inf_{
x\in\Omega_c,\ \mathbf u\in\mathbb U^N, 
\|\mathbf u-\mathbf u^\star(x)\|\ge R
}
\left\{
J_N(x,\mathbf u)-V_N(x)
\right\}.
\label{eq:uniform_separation_gap}
\end{equation}
For any $r>0$, define the admissible neighborhood of
the global minimizer and its minimum prior probability as
\begin{equation}
\mathcal B_{N,r}^{\mathrm{adm}}(x)
:=
B_r\!\left(\mathbf u^\star(x)\right)
\cap\mathbb U^N, \ 
p_{N,c}^{\mathrm{adm}}(r)
:=
\inf_{x\in\Omega_c}
P\!\left(\mathcal B_{N,r}^{\mathrm{adm}}(x)\right).
\label{eq:uniform_admissible_prior_coverage} 
\end{equation}
Here, $P(\mathcal B_{N,r}^{\mathrm{adm}}(x))$
denotes the probability that a control sequence
sampled from the prior $P$ lies within distance
$r$ of $\mathbf u^\star(x)$, while the infimum
ensures a uniform lower bound over $x\in\Omega_c$.

\begin{assumption}[Unique global minimizer]
\label{ass:unique_isolated_minimizer}
For every $c>0$ and every $x\in\Omega_c$,
the finite-horizon cost $J_N(x,\cdot)$ admits
a unique global minimizer $\mathbf u^\star(x)$.
Moreover, $\gamma_{N,c}(R)>0$ for every $R>0$.
\end{assumption}

\begin{assumption}[Uniform prior coverage]
\label{ass:uniform_prior_coverage}
For every $c>0$ and $r>0$,
we have $p_{N,c}^{\mathrm{adm}}(r)>0$.
\end{assumption}

\begin{remark}[Discussion of Assumptions 4 and 5]
Assumption~\ref{ass:unique_isolated_minimizer} ensures that control sequences sufficiently far from the global minimizer induces a uniformly positive finite-horizon cost gap $J_N(x,u)-V_N(x)$. Assumption~\ref{ass:uniform_prior_coverage} ensures that the prior assigns nonvanishing probability mass to the admissible neighborhoods of the global minimizer $\mathbf u^*(x)$, and is satisfied when $P$ has full support on the compact set $\mathbb U^N$, e.g. uniform distribution over $\mathbb U^N$. Together, these assumptions ensure that the Gibbs posterior increasingly concentrates around the global minimizer as $\lambda\to 0$.
\end{remark}

\begin{proposition}[Suboptimality of infinite-sample MPPI]
\label{prop:suboptimality_ideal_mppi}
Suppose Assumptions~\ref{ass:local_lipschitz_cost},
\ref{ass:unique_isolated_minimizer}, and
\ref{ass:uniform_prior_coverage} hold.
For any $N \in \mathbb N, c>0$, $\lambda>0$, $R>0$ and $r\in(0,R)$, define
\begin{equation}
\begin{aligned}
\varepsilon_{N,c}(\lambda)
:=
\inf_{\substack{R>0, r\in(0,R)\\g_{N,c}(r,R)>0}}
\Bigg\{
R+
\frac{2\sqrt{mN}u_{\max}}{
p_{N,c}^{\mathrm{adm}}(r)} \,
\exp\!\left(
-\frac{g_{N,c}(r,R)}{\lambda}
\right)
\Bigg\},
\end{aligned}
\label{eq:ideal_mppi_error_bound}
\end{equation}
where
$g_{N,c}(r,R):=\gamma_{N,c}(R)-L_{c,N}^u r$ with $L_{c,N}^u$ defined in Assumption \ref{ass:local_lipschitz_cost}.
Then, for all $x\in\Omega_c$,
\begin{equation}
\left\|
\mathbf u_{\infty,\lambda}(x)
-\mathbf u^\star(x)
\right\|
\leq \varepsilon_{N,c}(\lambda),
\label{eq:ideal_mppi_distance_to_global_optimum}
\end{equation}
\begin{equation}
J_N(x,\mathbf u_{\infty,\lambda}(x))-V_N(x)
\leq L_{c,N}^u\varepsilon_{N,c}(\lambda).
\label{eq:nonconvex_ideal_mppi_suboptimality}
\end{equation}
Moreover,
$\lim_{\lambda\to0}\varepsilon_{N,c}(\lambda)=0$.
\end{proposition}

\begin{remark}[Interpretation and comparison]
The radii $r$ and $R$ are auxiliary parameters introduced solely for the analysis and do not affect the MPPI algorithm. Since the bounds hold for any $0<r<R$ satisfying $g_{N,c}(r,R)>0$, taking the infimum yields the tightest bound within this family. For any fixed $R>0$, a sufficiently small $r$ ensures $g_{N,c}(r,R)>0$. The condition $g_{N,c}(r,R)>0$ ensures a positive cost separation between control sequences within distance $r$ of the global minimizer and those at least distance $R$ away. This separation ensures that the Gibbs weights of distant control sequences become exponentially smaller relative to those near the minimizer as $\lambda\to0$. As $\lambda\to0$, the auxiliary radii $R$ and $r$ can be chosen progressively smaller, while the decrease in temperature compensates for the shrinking separation gap and prior coverage, ensuring that $\varepsilon_{N,c}(\lambda)\to0$.

Unlike the contraction analysis in \citep{yi2024covo}, which relies on quadratic or additional smoothness and structural conditions, our result establishes convergence for general nonconvex and nonsmooth costs, without requiring strong convexity or differentiability. It also accommodates general sampling priors satisfying the uniform coverage condition.
\end{remark}

\subsection{Sample Complexity and Finite-Sample Approximation Error Bound} \label{subsec:sample_complexity_mppi}
We next quantify the finite-sample approximation error between the finite-sample MPPI solution \(\mathbf u_{M,\lambda}\) in \eqref{eq:finite_sample_mppi} and its infinite-sample counterpart \(\mathbf u_{\infty,\lambda}\) in \eqref{eq:ideal_mppi_control}. We establish a high-probability error bound that
characterizes the dependence on the sample size $M$,
temperature $\lambda$, cost variability and the discrepancy between
the Gibbs posterior $Q_{\lambda,x}^\star$ in
\eqref{eq:ideal_mppi_posterior} and the sampling prior $P$.

\begin{proposition}[Finite-sample MPPI approximation]
\label{prop:mppi_control_mismatch_hp}
For any $N\in\mathbb N$ and $x\in\mathbb R^n$, suppose Assumption~\ref{ass:convex_bounded_cost} holds. For any $\delta\in(0,1)$, the finite-sample MPPI control law generated by Algorithm~\ref{alg:mppi} with $M$ i.i.d. samples satisfies, with probability
at least $1-\delta$,
\begin{equation}
\begin{aligned}
\left\|
\mathbf u_{M,\lambda}(x)
-\mathbf u_{\infty,\lambda}(x)
\right\|
&\leq
4u_{\max}
\sqrt{\frac{mN}{M\delta}}\,
\exp\left(\frac{\Gamma_\lambda(x)}{2}\right) \leq
4u_{\max}
\sqrt{\frac{mN}{M\delta}}\,
\exp\left(\frac{C_N(x)^2}{2\lambda^2}\right),
\end{aligned}
\label{eq:mppi_control_mismatch_hp}
\end{equation}
\begin{equation}
\Gamma_\lambda(x)
:=
\log\left(1+\chi^2(Q_{\lambda,x}^\star\|P)\right) =
\frac{1}{\lambda^2}
\int_0^1\int_0^1
\operatorname{Var}_{Q_{\lambda, x,r+q}}
\left(J_N(x,\mathbf u)\right)
\,dr\,dq,
\label{eq:Gamma_t}
\end{equation}
with $Q_{\lambda,x}^\star$ being the Gibbs posterior \eqref{eq:ideal_mppi_posterior}, $C_N(x)$ being defined in Assumption \ref{ass:convex_bounded_cost} and $Q_{\lambda,x,\tau}(d\mathbf u) \propto e^{-\tau J_N(x,\mathbf u)/\lambda} P(d\mathbf u)$ for $\tau\in[0,2]$.
\end{proposition}

\begin{remark}[Interpretation of the finite-sample bound]
Compared with \citep{yoon2022sampling}, the bound
\eqref{eq:mppi_control_mismatch_hp} exhibits a milder
dimension dependence of $\sqrt{mN}$ rather than $mN$,
with an explicit $M^{-1/2}$ convergence rate.
The factor $\exp(\Gamma_\lambda(x)/2)$ captures the sampling
complexity through the discrepancy between the prior
$P$ and Gibbs posterior $Q_{\lambda,x}^\star$, equivalently
characterized by the cost variance along the Gibbs path as in \eqref{eq:Gamma_t}.
The bound \eqref{eq:mppi_control_mismatch_hp} reveals a fundamental trade-off: (i) smaller $\lambda$ complicates Monte Carlo estimation;
(ii) better alignment between the prior $P$ and Gibbs posterior
$Q_{\lambda,x}^*$ reduces approximation error; and
(iii) greater cost variability along the Gibbs path $Q_{\lambda,x,\tau}(d\mathbf u)$ increases
$\Gamma_\lambda(x)$ and the required sample size.
\end{remark}

\subsection{Semi-global Practical Recurrence} \label{subsec:stability}
This subsection establishes practical recurrence of finite-sample
MPPI by combining the suboptimality of infinite-sample MPPI
$\varepsilon_{N,c}(\lambda)$ in
\eqref{eq:ideal_mppi_error_bound} with the finite-sample
approximation error bound
\eqref{eq:mppi_control_mismatch_hp}.
We first derive a Lyapunov decrease condition for
the finite-horizon value function $V_N(x)$ defined in
\eqref{eq:optimal_value_function}, accounting for the
residual errors induced by the temperature $\lambda$
and finite sample size $M$. We then use this decrease
condition to establish semi-global practical recurrence \citep{GRAMMATICO20132939,moldenhauer2025robust} of the
closed-loop system under the finite-sample feedback law
$\kappa_{M,\lambda}$ defined in \eqref{eq:finite_sample_mppi_feedback}.

We first show the following single step Lyapunov decrease condition with high probability.
\begin{proposition}[High-probability Lyapunov decrease
of finite-sample MPPI]
\label{prop:mppi_finite_sample_drift}
Fix a prediction horizon $N\geq2$ and let
$c>0$. Consider the compact sublevel set
\(
\Omega_c:=\{x\in\mathbb R^n:V_N(x)\leq c\}.
\)
Suppose Assumptions~\ref{ass:convex_bounded_cost}--
\ref{ass:uniform_prior_coverage} hold. Then, for any $x\in\Omega_c, \tau>0$,
$\lambda>0$, $M\in\mathbb N$ and
$\delta\in(0,1)$ satisfying
\(
\Delta_{N,c}^{M}(x,\lambda,\delta)\leq\tau,
\)
the finite-sample MPPI control law satisfies
\begin{equation}
\begin{aligned}
\mathbb P\Big(
V_N(f(x,\kappa_{M,\lambda}(x)))-V_N(x) \leq-\alpha_{N,c,\tau}\sigma(x)
+\gamma_{c,\tau}
\Delta_{N,c}^{M}(x,\lambda,\delta)
\mid x
\Big)\geq1-\delta,
\end{aligned}
\label{eq:mppi_hp_drift}
\end{equation}
where, for $
\gamma_{c,\tau}
:=
\sup_{\sigma(x)\leq c+\tau}\bar a(x)<\infty$ and $\eta_{M,\lambda,N}(x,\delta)
:=
4u_{\max}
\sqrt{\frac{mN}{M\delta}}
\exp\left(\frac{C_N(x)^2}{2\lambda^2}\right)$,
\begin{equation}
\alpha_{N,c,\tau}
:=
1-(\gamma_{c,\tau}-1)^2
\left(
\frac{\gamma_{c,\tau}-1}{\gamma_{c,\tau}}
\right)^{N-2},
\label{eq:alpha_NC}
\end{equation}
\begin{equation}
\begin{aligned}
\Delta_{N,c}^{M}(x,\lambda,\delta)
:=L_{c,N}^u\Big[
&\varepsilon_{N,c}(\lambda)
+\eta_{M,\lambda,N}(x,\delta)
\Big],
\end{aligned}
\label{eq:Delta_finite_sample}
\end{equation}
with $\varepsilon_{N,c}(\lambda)$ defined in \eqref{eq:ideal_mppi_error_bound}, $L_{c,N}^u$ in Assumption~\ref{ass:local_lipschitz_cost} and
$\bar a(x)$ in Assumption~\ref{ass:no_terminal_cost_condition}.
\end{proposition}

For any \(\varepsilon,\tau>0\) and $N \in \mathbb N$, define the practical recurrence set
\(
\mathcal A_{c,\varepsilon}
\!:=\!
\left\{
x\in\Omega_c\!:\!
\alpha_{N,c,\tau}\sigma(x)
\leq
2\varepsilon
\right\}.
\)
We now establish practical recurrence of finite-sample MPPI by combining Proposition~\ref{prop:mppi_finite_sample_drift}.

\begin{proposition}[Practical recurrence of finite-sample MPPI]
\label{thm:mppi_semi-global_recurrence} 
Consider system~(1) controlled by
Algorithm~\ref{alg:mppi}.
Suppose Assumptions~1--5 hold. Then for any $r,\varepsilon>0$ and $p\in(0,1)$, 
\begin{equation}
\begin{gathered}
\exists\,T\in\mathbb N,\ c>0,\ 
\underline N\in\mathbb N,\,
\forall N\geq\underline N,\\
\exists\,\bar\lambda>0,\,
\forall\lambda\in(0,\bar\lambda),\quad
\exists\,\underline M\in\mathbb N,\,
\forall M\geq\underline M:\\
\mathbb P\left[
\exists\,k\in\{0,\ldots,T\}:
x_k\in\mathcal A_{c,\varepsilon}
\right]>1-p, \forall\,x_0\in\mathbb R^n,\, \|x_0\|\leq r.
\end{gathered}
\label{eq:mppi_finite_time_recurrence}
\end{equation}
\end{proposition}

\begin{remark}[Interpretation of practical recurrence]
Proposition~\ref{thm:mppi_semi-global_recurrence} establishes semi-global
practical recurrence in the sense of \citep{GRAMMATICO20132939,moldenhauer2025robust}.
Specifically, for any prescribed bounded set of initial conditions and any desired probability level, the closed-loop trajectory enters the practical set $\mathcal A_{c,\varepsilon}$ within a finite time $T$ with sufficiently high probability $1-p$, provided that the prediction horizon $N$ and sample size $M$ are sufficiently large and the temperature $\lambda$ is sufficiently small. The term ``semi-global'' refers to the arbitrary choice of a bounded initial set of size $r$, while ``practical'' refers to recurrence to a neighborhood of size $\epsilon$ rather than the origin.

Proposition~\ref{thm:mppi_semi-global_recurrence} provides a natural and mathematically appropriate closed-loop guarantee in the presence of persistent stochastic approximation errors. Finite-sample MPPI introduces fresh randomness at every control update, as characterized by Proposition~\ref{prop:mppi_control_mismatch_hp}. Consequently, even when the approximation error is small with high probability at each step, occasional inaccurate control updates cannot generally be excluded over an infinite time horizon. This prevents us from directly extending a high-probability one-step Lyapunov decrease condition to an infinite-time stability guarantee. Practical recurrence instead captures the achievable closed-loop behavior under such persistent uncertainty: for any prescribed finite-time confidence level, the trajectory remains bounded and reaches a neighborhood of the equilibrium within a finite time with arbitrarily high probability. Moreover, both the practical neighborhood $\epsilon$ and the failure probability $p$ can be made arbitrarily small through appropriate choices of the prediction horizon $N$, temperature $\lambda$, and sample size $M$. Thus, rather than requiring stochastic approximation errors to vanish permanently, the result establishes a tunable probabilistic guarantee on the closed-loop behavior, consistent with the motivation for practical recurrence under persistent stochastic perturbations in \citep{GRAMMATICO20132939,moldenhauer2025robust}.
\end{remark}

\section{Simulation}

\begin{figure*}[!t]
\centering
\begin{minipage}[t]{0.48\textwidth}
\centering
\includegraphics[width=.48\linewidth]{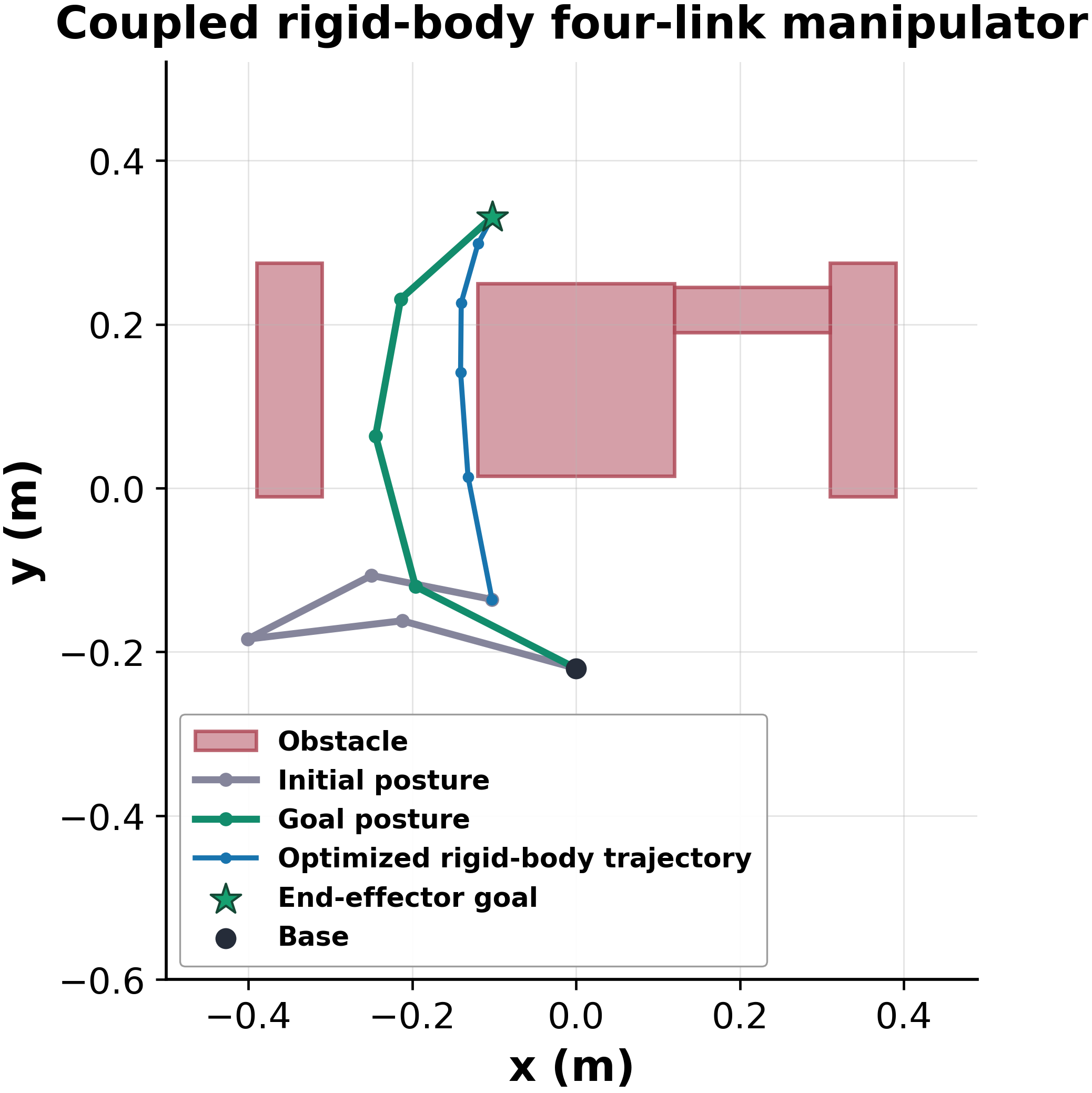}
\hfill
\includegraphics[width=.48\linewidth]{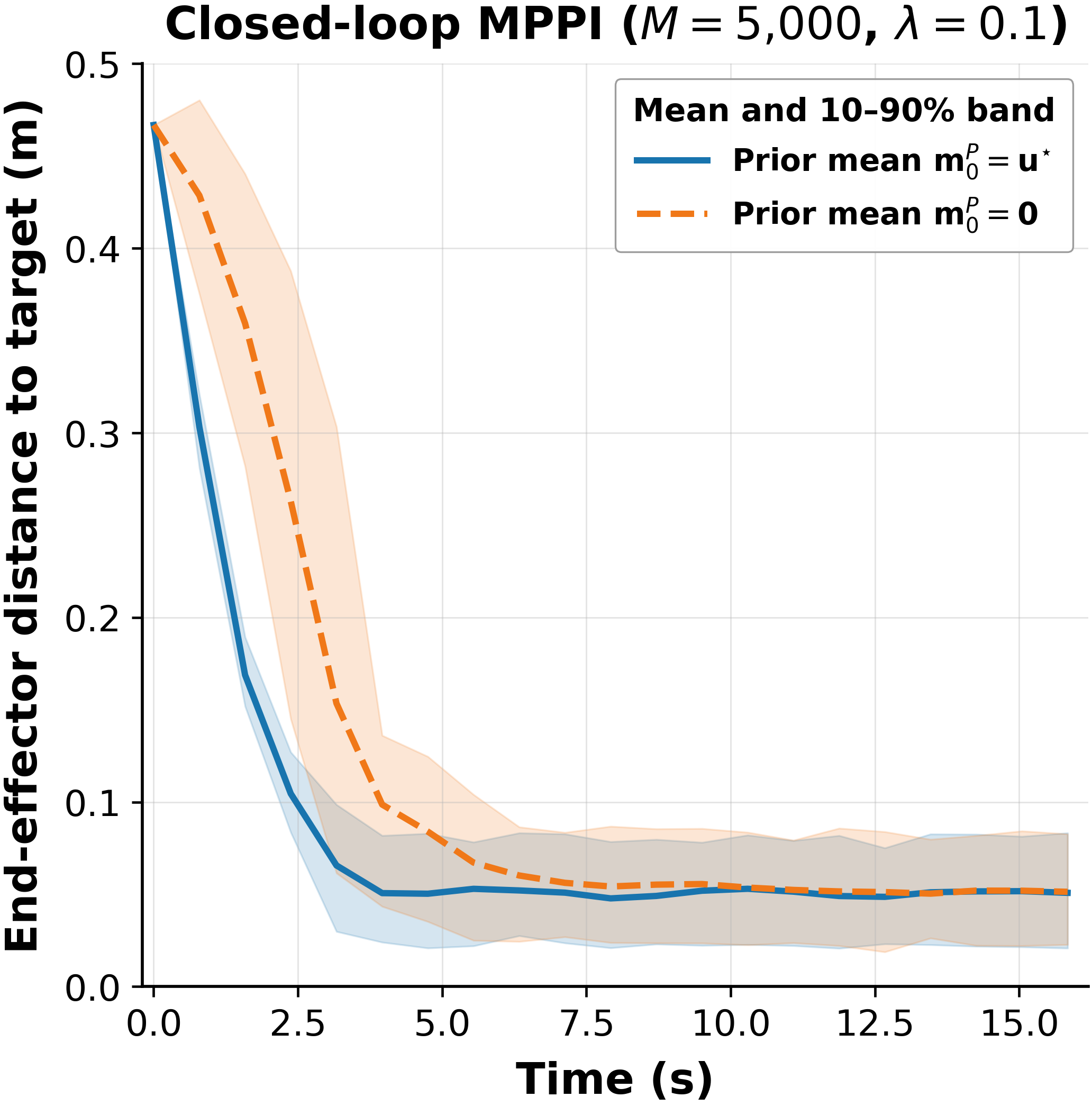}
\captionof{figure}{Environment and closed-loop trajectories.}
\label{fig:env_setup}
\end{minipage}
\hfill
\begin{minipage}[t]{0.48\textwidth}
\centering
\includegraphics[width=.48\linewidth]{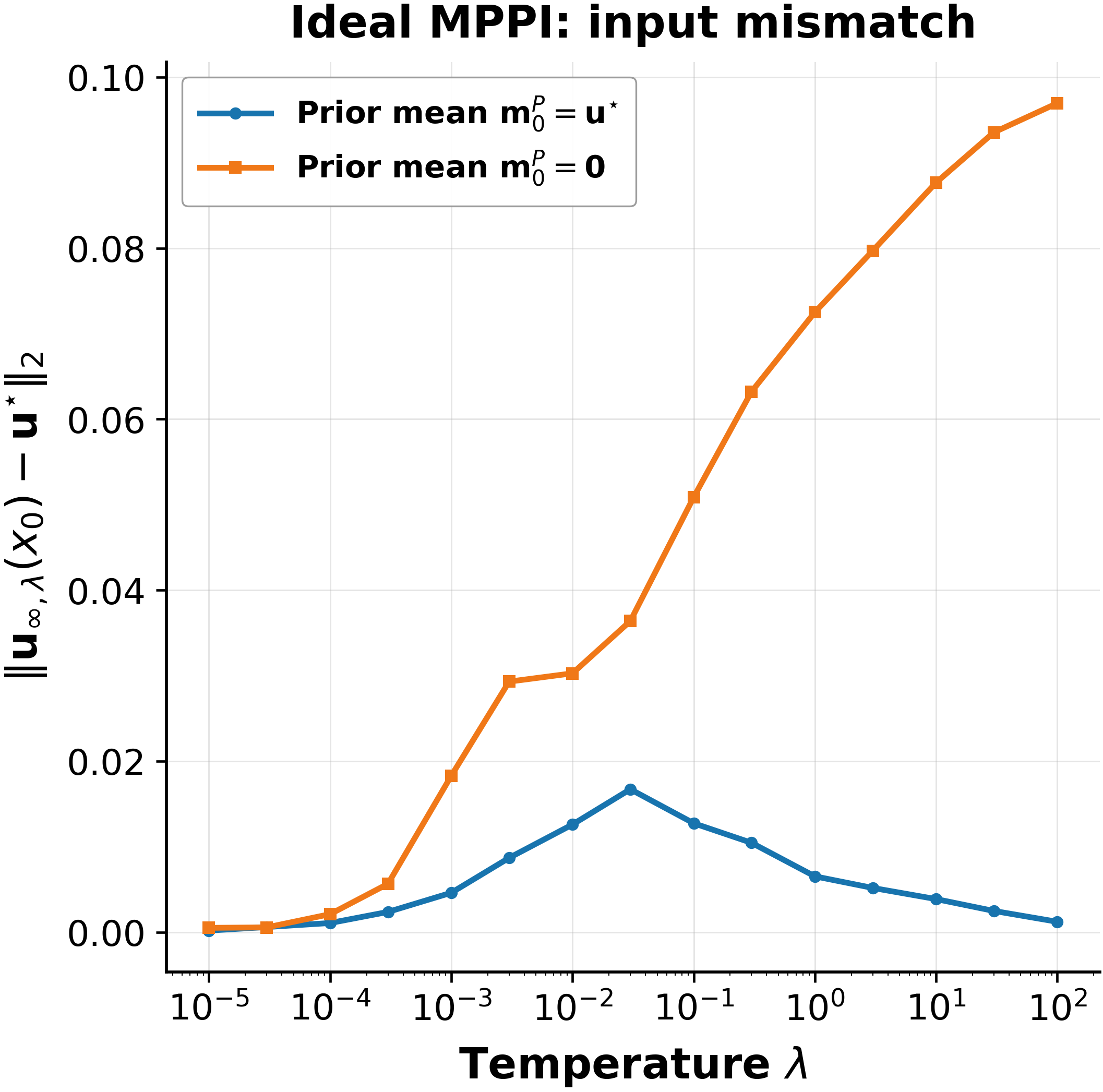}
\hfill
\includegraphics[width=.48\linewidth]{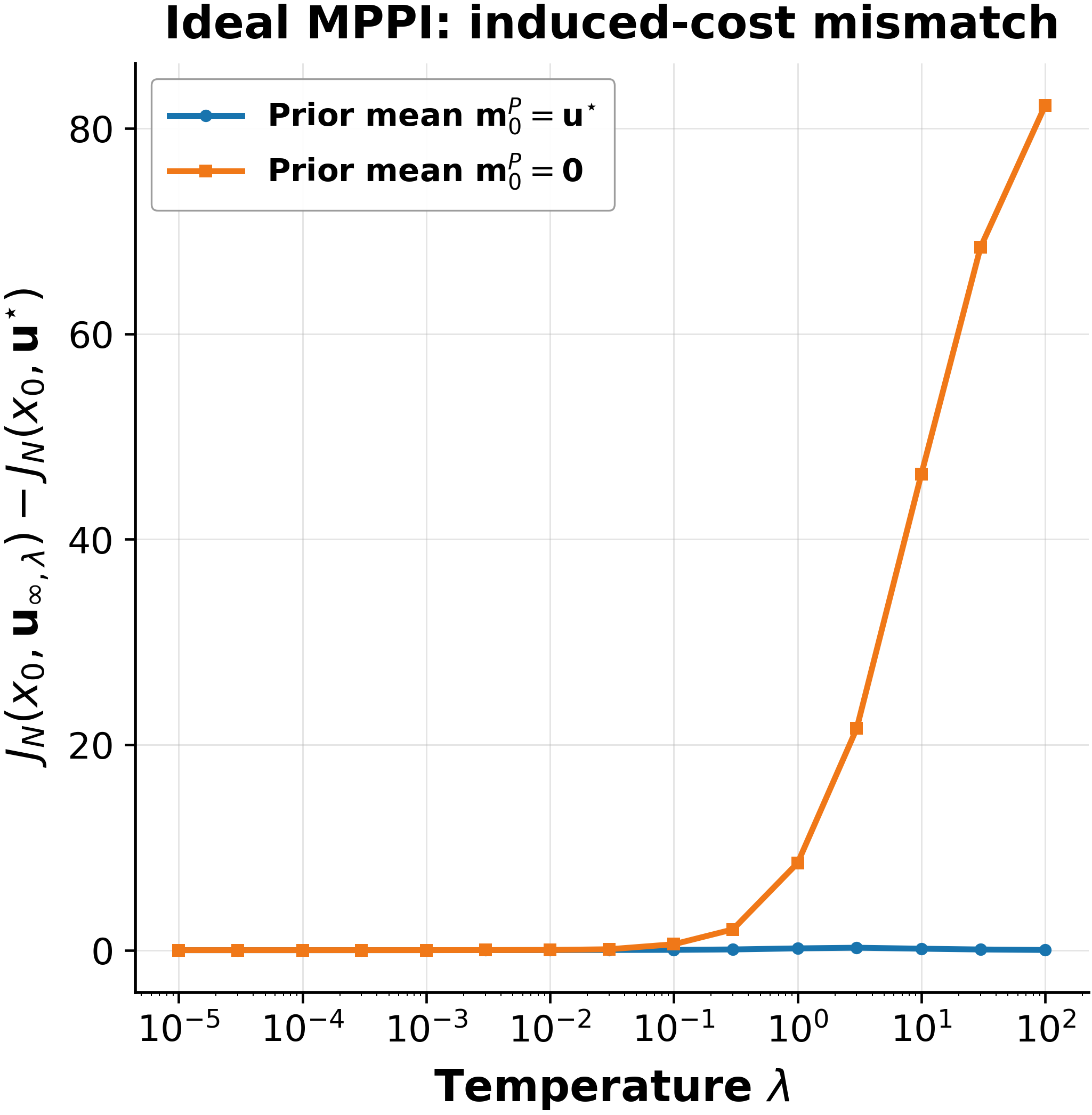}
\captionof{figure}{Temperature-dependent ideal-MPPI suboptimality.}
\label{fig:lbd_ablation}
\end{minipage}

\vspace{0.6em}

\begin{minipage}[t]{0.48\textwidth}
\centering
\includegraphics[width=.48\linewidth]{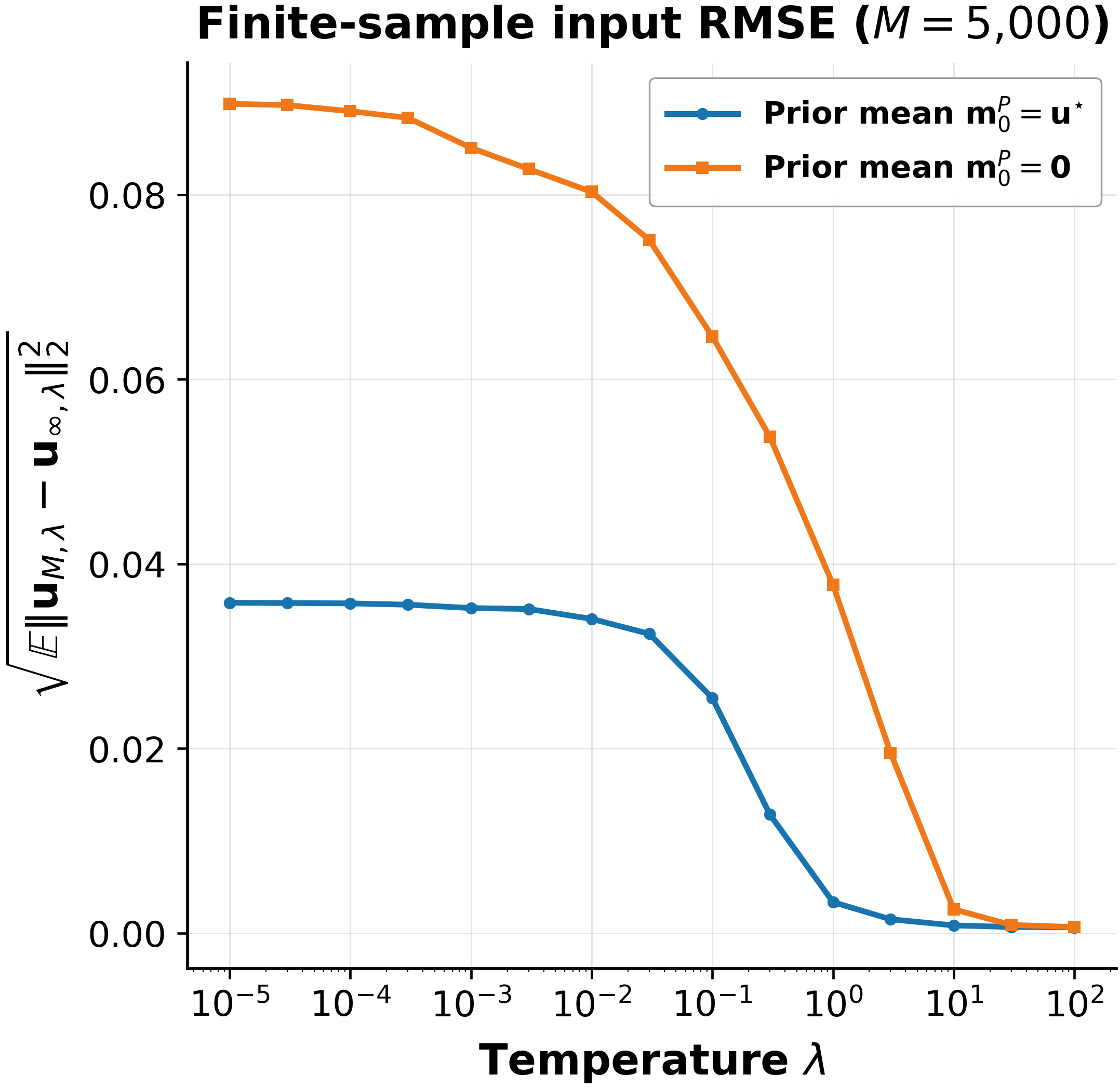}
\hfill
\includegraphics[width=.48\linewidth]{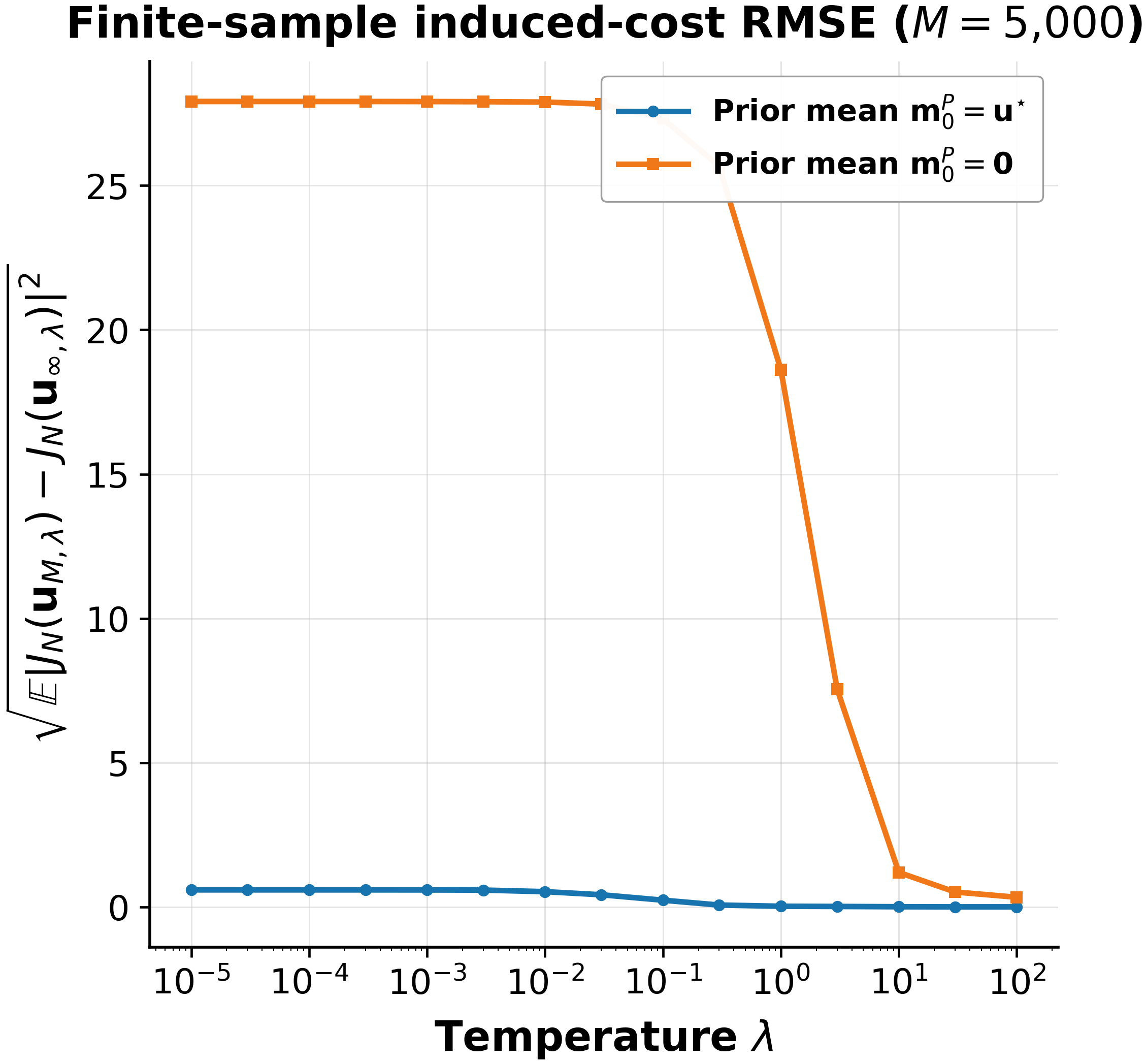}
\captionof{figure}{Temperature-dependent finite-sample error.}
\label{fig:lbd_ablation_finite_sample}
\end{minipage}
\hfill
\begin{minipage}[t]{0.48\textwidth}
\centering
\includegraphics[width=.48\linewidth]{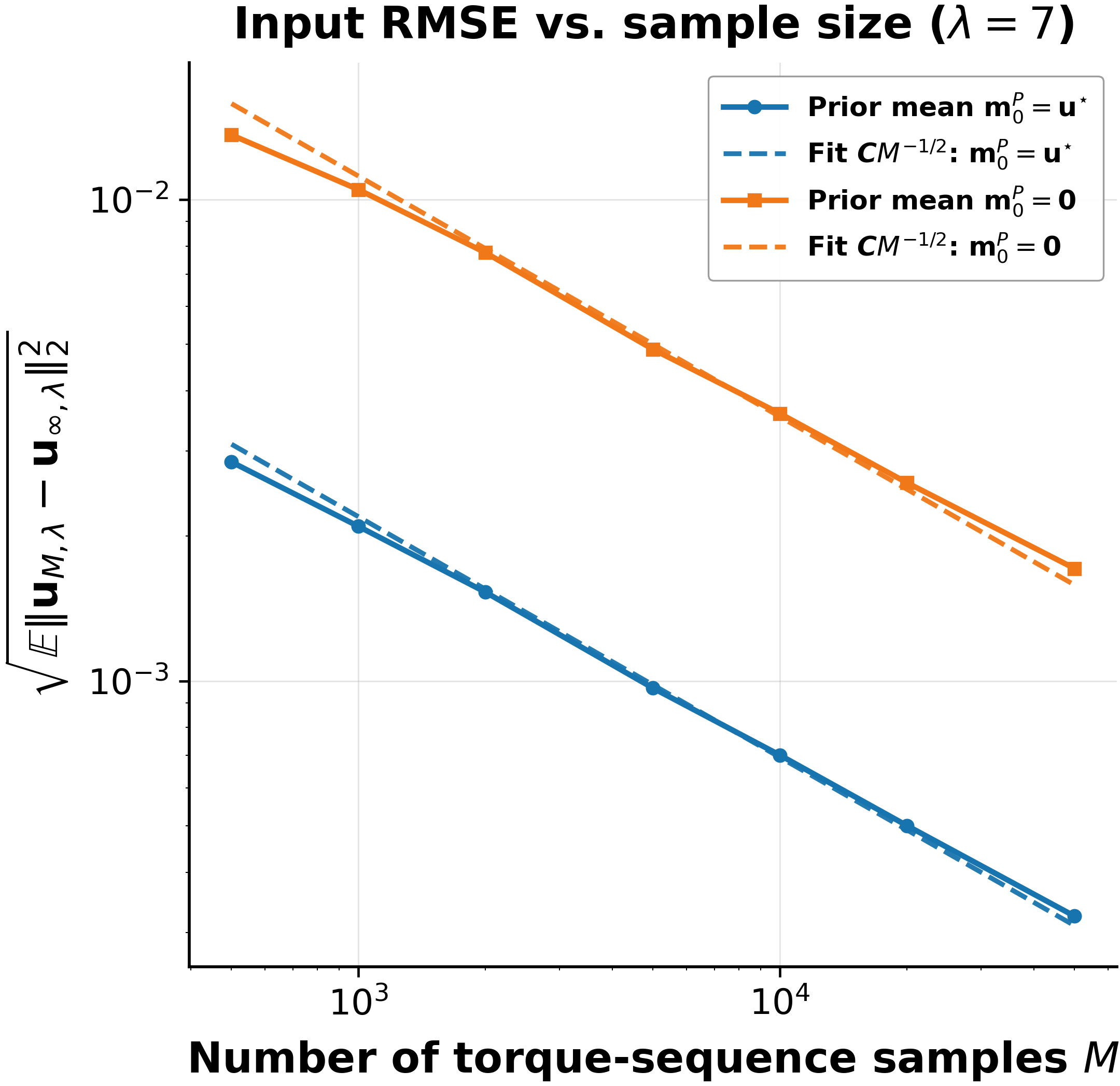}
\hfill
\includegraphics[width=.48\linewidth]{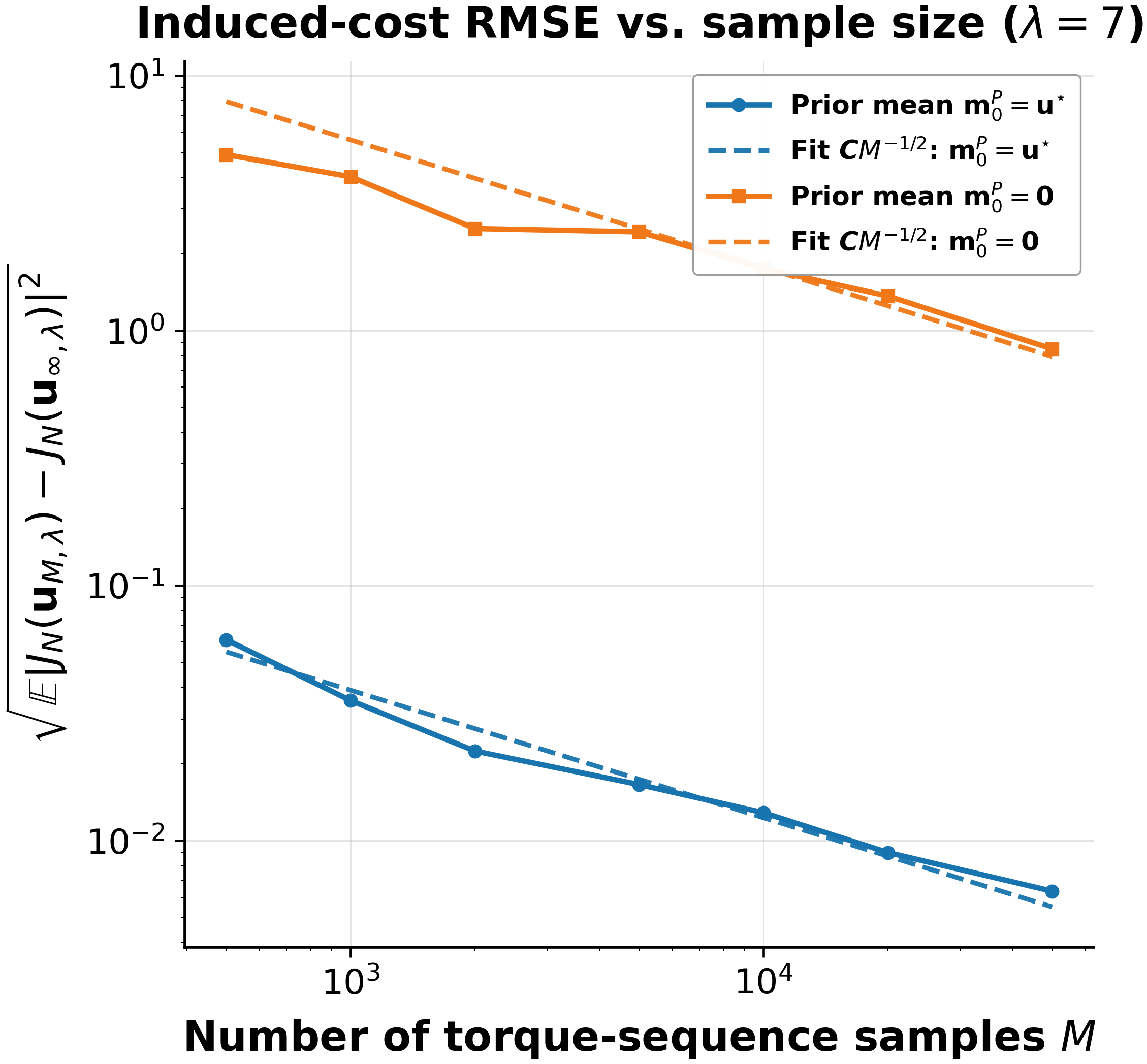}
\captionof{figure}{Finite-sample convergence with $M$.}
\label{fig:M_ablation}
\end{minipage}
\end{figure*}

We consider a four-link planar manipulator navigating a narrow
passage with obstacles (Fig.~\ref{fig:env_setup}) \citep{colbaugh1989cartesian,colbaugh1989obstacle}, governed by
\(
    \dot q=v,
    M(q)\dot v+C(q,v)v+Dv=u,
\)
where $q,v,u\in\mathbb R^4$ denote joint angles, angular velocities,
and torques. The links are uniform rods with lengths
$\ell=(0.22,0.19,0.17,0.15)^\top$ m, linear density
$1$ kg/m, and viscous damping
$D=\operatorname{diag}(0.060,0.050,0.040,0.030)$.
The initial and target configurations are
\(
q_0=(2.8745,0.3852,-2.7850,-0.6711)^\top, v_0=0,
q_{\mathrm{goal}}=(2.67,-0.84,-0.44,-0.66)^\top,
p_{\mathrm{goal}}=(-0.1023,0.3309)^\top\ \mathrm{m}.
\)
We use $N=5$, $\Delta t=0.792$ s, and semi-implicit Euler
integration with 24 substeps per interval. The input bounds are
$|u_k|\leq(0.10,0.085,0.075,0.060)^\top$ componentwise.
The finite-horizon objective is
\[
\begin{aligned}
J_N(x_0,\mathbf u)
={}&\sum_{k=0}^{N-1}\Delta t
\bigl(
0.2\|p(q_{k+1})-p_{\mathrm{goal}}\|^2
+0.025\|v_{k+1}\|^2
+0.006\|u_k\|^2
\bigr)\\
&+40\|p(q_N)-p_{\mathrm{goal}}\|^2
+12\|q_N-q_{\mathrm{goal}}\|^2
+\|v_N\|^2+\Phi(\mathbf u),
\end{aligned}
\]
where $\Phi$ penalizes obstacle proximity, self-collision,
and joint-limit violations. An empirical minimizer $\mathbf u^\star$ is obtained through offline multistart optimization and local refinement using L-BFGS-B solver. We consider truncated Gaussian priors
$P=\operatorname{TN}_{\mathbb U^N}(m_0^P,10^{-4}I_{20})$,
with well-aligned ($m_0^P=\mathbf u^\star$) and poorly
aligned ($m_0^P=0$) means. An empirical value of $\mathbf u_{\infty,\lambda}$ is estimated by using $400,000$ samples for given $\lambda>0$. Four sets of results are illustrated below. (1) Figure~\ref{fig:env_setup} shows closed-loop MPPI with
$M=5000$, $\lambda=0.1$, and 200 Monte Carlo trajectories
per prior over 20 control intervals. The tip-to-target
distance decreases and remains within a small residual
neighborhood, illustrating empirical practical convergence.
Shaded regions indicate the 10--90\% ensemble bands. (2) To illustrate Proposition~\ref{prop:suboptimality_ideal_mppi},
we vary $\lambda$ over 15 values from $10^{-5}$ to $100$,
keeping the prior covariance fixed. We evaluate
\(
e_u(\lambda)
=\|\mathbf u_{\infty,\lambda}-\mathbf u^\star\|,
e_J(\lambda)=J_N(x_0,\mathbf u_{\infty,\lambda})
-J_N(x_0,\mathbf u^\star).
\)
Figure~\ref{fig:lbd_ablation} shows smaller errors at low
temperatures and generally better performance with the
well-aligned prior. Its nonmonotonic input mismatch reflects
the competition between Gibbs concentration as
$\lambda\to0$ and convergence to the prior mean as
$\lambda\to\infty$. (3) We fix $M=5000$ and evaluate 15 temperatures using
20 independent Monte Carlo realizations per prior.
The finite-sample error is measured by
\(
e_{\mathrm{RMSE}}(\lambda)
=\sqrt{\mathbb E\left[
\|\mathbf u_{M,\lambda}
-\mathbf u_{\infty,\lambda}\|^2
\right]}
\)
, estimated using 20 independent Monte Carlo
realizations for each prior and temperature. Figure~\ref{fig:lbd_ablation_finite_sample} shows that
smaller $\lambda$ generally increases input and induced-cost
errors, consistent with the temperature dependence in
\eqref{eq:gibbs_variance_bound_boundedJ}
and \eqref{eq:mppi_control_mismatch_hp}. (4) We fix $\lambda=7$ and vary
$M\in\{500,1000,2000,5000,10000,20000,50000\}$,
using 24 Monte Carlo realizations per setting.
The numerical infinite-sample reference uses
$4\times10^5$ prior samples per prior.
Figure~\ref{fig:M_ablation} shows that both input and
induced-cost RMSE decrease approximately as $M^{-1/2}$,
consistent with \eqref{eq:mppi_control_mismatch_hp}.
The well-aligned prior consistently yields smaller errors.

\section{Conclusions}
This paper investigated the suboptimality, sample complexity and closed-loop guarantees of MPPI under finite sample and positive temperature, with suboptimality guarantees established under local Lipschitz continuity without requiring convexity or differentiability. The analysis reveals an important tradeoff: reducing the temperature improves approximation of the original optimum but increases the sample complexity. The sampling distribution and its coverage of near-optimal solutions therefore play a critical role in MPPI performance. By incorporating these approximation errors into the MPC value-function decrease, semi-global practical recurrence can be established under appropriate parameter selections. Numerical results illustrate these effects and the importance of prior quality. Future work will consider model uncertainty and stronger stability results.

\section*{Generative AI Use}

Generative AI tools were used to assist with mathematical derivations checking, experimental implementation, literature review, and manuscript editing. The authors reviewed and verified the AI-assisted content and take full responsibility
for the accuracy and integrity of the final manuscript.

\acks{
The work of Changyi Lei was supported by the Australian
Government Research Training Program (RTP) Scholarship
and the Australian Research Council under Project
DE220101527. The work of Dragan Ne\v{s}i\'{c} was
supported by the Australian Research Council through
Discovery Project DP250100300. The work of Ye Pu was
supported by the Australian Research Council under
Project DP260101027.
}

\bibliography{l4dc2027-sample}

\newpage
\appendix
\section*{Supporting Theoretical Results}
\begin{figure}[htbp]
    \centering
    \includegraphics[width=0.8\textwidth]{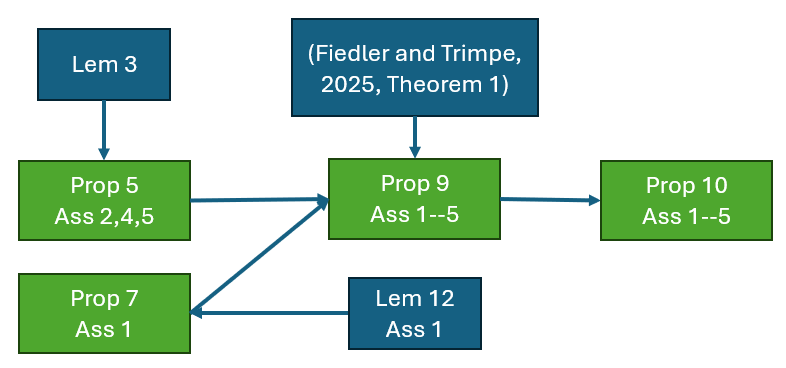}
    \caption{Theoretical results dependence graph.}
    \label{fig:example}
\end{figure}

This section lists some supportive results that are needed for the proof of the main results. Figure \ref{fig:example} shows the dependence relation between all theoretical results in this paper.

\begin{lemma}
    Using the same notation as in Proposition \ref{prop:mppi_control_mismatch_hp}, we have
    \begin{equation}
\begin{split}
    \operatorname{Var}\left[\overline r_M  \right] = \frac{\chi^2(Q_{\lambda,x}^* \| P)}{M},
\end{split} \label{eq:var_rM_chi_divergence}
\end{equation}
\begin{equation}
    \operatorname{Var}
    \left(
        \overline r_M
    \right)
    =
    \frac{
        \exp(\Gamma_\lambda(x))-1
    }{
        M
    }.
    \label{eq:var_rM_var_J}
\end{equation}  
\label{lem:var_rM_var_J}
\end{lemma}

\begin{proof}
    Equation \eqref{eq:var_rM_chi_divergence} follows directly from the definition of $\chi^2$ divergence and \eqref{eq:var_rM}.

We next express the $\chi^2$ divergence in terms of the variance of the
finite-horizon cost along a Gibbs path, which shows \eqref{eq:var_rM_var_J}. For $\tau\in[0,2]$, define
\begin{equation}
    Q_{\lambda,x,\tau}(d\mathbf u)
    :=
    \frac{
        \exp\!\left(
            -\tau J_N(x,\mathbf u)/\lambda
        \right)
    }{
        Z_{\lambda,x}(\tau)
    }
    P(d\mathbf u),
    \label{eq:gibbs_path_t}
\end{equation}
where
\begin{equation}
    Z_{\lambda,x}(\tau)
    :=
    \mathbb E_P
    \left[
        \exp\!\left(
            -\tau J_N(x,\mathbf u)/\lambda
        \right)
    \right].
\end{equation}
By the definition of the $\chi^2$ divergence,
\begin{align}
    1+\chi^2(Q_{\lambda,x}^*\|P)
    &=
    \mathbb E_P
    \left[
        \left(
            \frac{dQ_{\lambda,x}^*}{dP}(\mathbf u)
        \right)^2
    \right]
    \notag\\
    &=
    \frac{1}{Z_{\lambda,x}(1)^2}
    \mathbb E_P
    \left[
        \exp\!\left(
            -\frac{2J_N(x,\mathbf u)}{\lambda}
        \right)
    \right]
    \notag\\
    &=
    \frac{
        Z_{\lambda,x}(2)
    }{
        Z_{\lambda,x}(1)^2
    }.
    \label{eq:chi_partition_t}
\end{align}

Define the log-partition function
\begin{equation}
    \phi_{\lambda,x}(\tau)
    :=
    \log Z_{\lambda,x}(\tau),
\end{equation}
where
\begin{equation}
    Z_{\lambda,x}(\tau)
    =
    \mathbb E_P
    \left[
        \exp\left(
            -\frac{\tau}{\lambda}
            J_N(x,\mathbf u)
        \right)
    \right].
\end{equation}
We first compute the first derivative of $Z_{\lambda,x}(\tau)$. Differentiating
under the expectation gives
\begin{align}
    Z_{\lambda,x}'(\tau)
    &=
    \frac{d}{d\tau}
    \mathbb E_P
    \left[
        \exp\left(
            -\frac{\tau}{\lambda}
            J_N(x,\mathbf u)
        \right)
    \right]
    \nonumber\\
    &=
    \mathbb E_P
    \left[
        \frac{d}{d\tau}
        \exp\left(
            -\frac{\tau}{\lambda}
            J_N(x,\mathbf u)
        \right)
    \right]
    \nonumber\\
    &=
    -\frac{1}{\lambda}
    \mathbb E_P
    \left[
        J_N(x,\mathbf u)
        \exp\left(
            -\frac{\tau}{\lambda}
            J_N(x,\mathbf u)
        \right)
    \right].
\label{eq:Z_tau_first_derivative}
\end{align}
We have
\begin{align}
    \phi_{\lambda,x}'(\tau)
    &=
    \frac{Z_{\lambda,x}'(\tau)}
         {Z_{\lambda,x}(\tau)}
     =
    -\frac{1}{\lambda}
    \frac{
        \mathbb E_P
        \left[
            J_N(x,\mathbf u)
            \exp\left(
                -\frac{\tau}{\lambda}
                J_N(x,\mathbf u)
            \right)
        \right]
    }{
        Z_{\lambda,x}(\tau)
    }.
\label{eq:phi_tau_first_derivative_expanded}
\end{align}
By \eqref{eq:gibbs_path_t}, we have equivalently
\begin{equation}
    \phi_{\lambda,x}'(\tau)
    =
    -\frac{1}{\lambda}
    \mathbb E_{Q_{\lambda,x,\tau}}
    \left[
        J_N(x,\mathbf u)
    \right].
\label{eq:phi_tau_first_derivative}
\end{equation}

Next, differentiating \eqref{eq:Z_tau_first_derivative} once more gives
\begin{align}
    Z_{\lambda,x}''(\tau)
    &=
    -\frac{1}{\lambda}
    \frac{d}{d\tau}
    \mathbb E_P
    \left[
        J_N(x,\mathbf u)
        \exp\left(
            -\frac{\tau}{\lambda}
            J_N(x,\mathbf u)
        \right)
    \right]
    \nonumber\\
    &=
    -\frac{1}{\lambda}
    \mathbb E_P
    \left[
        J_N(x,\mathbf u)
        \frac{d}{d\tau}
        \exp\left(
            -\frac{\tau}{\lambda}
            J_N(x,\mathbf u)
        \right)
    \right]
    \nonumber\\
    &=
    \frac{1}{\lambda^2}
    \mathbb E_P
    \left[
        J_N(x,\mathbf u)^2
        \exp\left(
            -\frac{\tau}{\lambda}
            J_N(x,\mathbf u)
        \right)
    \right].
\label{eq:Z_tau_second_derivative}
\end{align}
The second derivative of $\phi_{\lambda,x}(\tau)$ is
\[
\begin{aligned}
    \phi_{\lambda,x}''(\tau)
    &=
    \frac{
        Z_{\lambda,x}''(\tau)
    }{
        Z_{\lambda,x}(\tau)
    }
    -
    \left(
        \frac{
            Z_{\lambda,x}'(\tau)
        }{
            Z_{\lambda,x}(\tau)
        }
    \right)^2
    \nonumber\\
    &=
    \frac{1}{\lambda^2}
    \frac{
        \mathbb E_P
        \left[
            J_N(x,\mathbf u)^2
            \exp\left(
                -\frac{\tau}{\lambda}
                J_N(x,\mathbf u)
            \right)
        \right]
    }{
        Z_{\lambda,x}(\tau)
    }
    -
    \frac{1}{\lambda^2}
    \left(
        \frac{
            \mathbb E_P
            \left[
                J_N(x,\mathbf u)
                \exp\left(
                    -\frac{\tau}{\lambda}
                    J_N(x,\mathbf u)
                \right)
            \right]
        }{
            Z_{\lambda,x}(\tau)
        }
    \right)^2.
\end{aligned}
\]
Again using the definition of $Q_{\lambda,x,\tau}$, we obtain
\begin{align}
    \phi_{\lambda,x}''(\tau)
    &=
    \frac{1}{\lambda^2}
    \mathbb E_{Q_{\lambda,x,\tau}}
    \left[
        J_N(x,\mathbf u)^2
    \right]
    -
    \frac{1}{\lambda^2}
    \left(
        \mathbb E_{Q_{\lambda,x,\tau}}
        \left[
            J_N(x,\mathbf u)
        \right]
    \right)^2
    \nonumber\\
    &=
    \frac{1}{\lambda^2}
    \left\{
        \mathbb E_{Q_{\lambda,x,\tau}}
        \left[
            J_N(x,\mathbf u)^2
        \right]
        -
        \left(
            \mathbb E_{Q_{\lambda,x,\tau}}
            \left[
                J_N(x,\mathbf u)
            \right]
        \right)^2
    \right\}
    \nonumber\\
    &=
    \frac{1}{\lambda^2}
    \operatorname{Var}_{Q_{\lambda,x,\tau}}
    \left(
        J_N(x,\mathbf u)
    \right).
\end{align}
Therefore,
\begin{equation}
    \phi_{\lambda,x}''(\tau)
    =
    \frac{1}{\lambda^2}
    \operatorname{Var}_{Q_{\lambda,x,\tau}}
    \left(
        J_N(x,\mathbf u)
    \right).
\label{eq:phi_tau_second_derivative}
\end{equation}

Taking logarithms in \eqref{eq:chi_partition_t},
\begin{align}
    \log &
    \left(
        1+\chi^2(Q_{\lambda,x}^*\|P)
    \right)
    =
    \phi_{\lambda,x}(2)
    -
    2\phi_{\lambda,x}(1)
    +
    \phi_{\lambda,x}(0)
    \notag\\
    &=
    \int_0^1
    \left[
        \phi_{\lambda,x}'(1+q)
        -
        \phi_{\lambda,x}'(q)
    \right]
    dq
    \notag\\
    &=
    \int_0^1
    \int_0^1
    \phi_{\lambda,x}''(r+q)
    \,dr\,dq
    \notag\\
    &=
    \frac{1}{\lambda^2}
    \int_0^1
    \int_0^1
    \operatorname{Var}_{Q_{\lambda, x,r+q}}
    \left(
        J_N(x,\mathbf u)
    \right)
    \,dr\,dq.
\end{align}
Therefore,
\begin{equation}
    1+\chi^2(Q_{\lambda,x}^*\|P)
    =
    \exp(\Gamma_\lambda(x)).
\end{equation}
Combining with \eqref{eq:var_rM_chi_divergence} renders
\begin{equation}
    \operatorname{Var}
    \left(
        \overline r_M
    \right)
    =
    \frac{
        \exp(\Gamma_\lambda(x))-1
    }{
        M
    }.
\end{equation}
\end{proof}

\section*{Proof of Lemma \ref{lem:ideal_mppi_variational}}
Consider the KL-regularized objective
\[
\mathbb E_{\mathbf u\sim Q}[J_N(x,\mathbf u)]
+\lambda\,\mathrm{KL}(Q\|P).
\]
By the Gibbs variational principle
\cite[Lemma 2.2]{alquier2024user}, its unique minimizer
$Q_{\lambda,x}^\star$ is given by \eqref{eq:ideal_mppi_posterior}.
Consequently,
\begin{align*}
\mathbb E_{\mathbf u\sim Q_{\lambda,x}^\star}[\mathbf u]
&=
\frac{
\mathbb E_{\mathbf u\sim P}
\left[e^{-J_N(x,\mathbf u)/\lambda}\mathbf u\right]
}{
\mathbb E_{\mathbf u\sim P}
\left[e^{-J_N(x,\mathbf u)/\lambda}\right]
} =\mathbf u_{\infty,\lambda}(x),
\end{align*}
where the last equality follows directly from
\eqref{eq:ideal_mppi_control}.

\section*{Proof of Proposition \ref{prop:suboptimality_ideal_mppi}}
    The proof steps are inspired by \citep[Proposition 4.5]{fornasier2024consensus}. The intuition of the proof is that if the prior puts nontrivial mass near the global minimizer $\mathbf u^\star(x)$, then so long as $\lambda$ is small enough, the MPPI algorithm will prefer those inputs around $\mathbf u^\star(x)$ and so is the final $\mathbf u_{\infty,\lambda}$ and the induced suboptimality gap in $J_N(x, \cdot)$.

    \textbf{(1) } For simplicity, denote $\mathbf u^\star:= \mathbf u^\star(x)$ and define, for any $a \geq 0$, \(
\mathcal B_{N,a}^{\mathrm{adm}}(x)
:=
B_a\!\left(\mathbf u^\star\right)
\cap
\mathbb U^N. \)
Assumption \ref{ass:local_lipschitz_cost} implies that
\begin{equation}
J_N(x,\mathbf u)
\leq
V_N(x)+L_{c,N}^u r, \forall \,
\mathbf u
\in
\mathcal B_{N,r}^{\mathrm{adm}}(x).
\label{eq:local_upper_bound_near_optimum}
\end{equation}
Recall the posterior definition \eqref{eq:ideal_mppi_posterior} and denote $Z_{\lambda,x}:= \mathbb E_{\mathbf u\sim P}\!
\left[
\exp\!\left(-J_N(x,\mathbf u)/\lambda\right)
\right]$, then 
\begin{align}
Z_{\lambda,x}
& \geq
\int_{\mathcal B_{N,r}^{\mathrm{adm}}}
\exp\!\left(
-\frac{J_N(x,\mathbf u)}{\lambda}
\right)
P(d\mathbf u)
\overset{\eqref{eq:local_upper_bound_near_optimum}}{\geq}
P\!\left(
\mathcal B_{N,r}^{\mathrm{adm}}
\right)
\exp\!\left(
-\frac{V_N(x)+L_{c,N}^u r}{\lambda}
\right)
\nonumber\\
&\overset{\eqref{eq:uniform_admissible_prior_coverage}}{\geq}
p_{N,c}^{\mathrm{adm}}(r)
\exp\!\left(
-\frac{V_N(x)+L_{c,N}^u r}{\lambda}
\right).
\label{eq:gibbs_partition_lower_bound}
\end{align}
where the first inequality comes from Assumption \ref{ass:uniform_prior_coverage} that $p_{N,c}^{\mathrm{adm}}(r)>0.$ On the other hand, by the definition of
$\gamma_{N,c}(R)$ as in \eqref{eq:uniform_separation_gap} and Assumption \ref{ass:unique_isolated_minimizer}, every admissible sequence satisfying
\(
\|\mathbf u-\mathbf u^\star\|\geq R
\)
satisfies
\begin{equation}
J_N(x,\mathbf u)
\geq
V_N(x)+\gamma_{N,c}(R).
\label{eq:cost_lower_bound_away_optimum}
\end{equation}

Define the event
\(
    \mathcal E_R:= \{ \mathbf u \in \mathcal B_{N,R}^{\mathrm{adm}}(x)\}
\)
with its complement $\mathcal E_R^c$. From \eqref{eq:ideal_mppi_posterior} in Lemma \ref{lem:ideal_mppi_variational}, we have
\begin{align}
Q_{\lambda,x}^*\!\left(\mathcal E_R^c\right)  &=
\int_{
\mathcal E_R^c
}
\frac{\exp\!\left(
-J_N(x,\mathbf u)/ \lambda
\right)}{Z_{\lambda,x}}
P(d\mathbf u)
 \overset{\eqref{eq:cost_lower_bound_away_optimum}}{\leq}
\frac{\exp\!\left(
-(V_N(x)+\gamma_{N,c}(R))/ \lambda
\right)}{Z_{\lambda,x}}
\nonumber\\
&\overset{\eqref{eq:gibbs_partition_lower_bound}}{\leq}
\frac{\exp\!\left(
-V_N(x)+\gamma_{N,c}(R)/ \lambda
\right)}{p_{N,c}^{\mathrm{adm}}(r)
\exp\!\left(
-(V_N(x)+L_{c,N}^u r)/\lambda
\right)}  \leq \frac{1}{p_{N,c}^{\mathrm{adm}}(r)}\exp\!\left(
-\frac{g_{N,c}(r,R)}{\lambda}
\right).
\label{eq:gibbs_tail_numerator_bound}
\end{align}

\textbf{(2)} We now proceed to show \eqref{eq:ideal_mppi_distance_to_global_optimum}. We have
\begin{align}
\!\!\!\left\|
\mathbf u_{\infty,\lambda}
-
\mathbf u^\star
\right\|
&=
\left\|
\mathbb E_{Q_{\lambda,x}^*}\!
[
\mathbf u-\mathbf u^\star
]
\right\| \leq
\mathbb E_{Q_{\lambda,x}^*}\!\!
\left[
\left\|
\mathbf u-\mathbf u^\star
\right\|
\right] \nonumber\\
&=\mathbb E_{Q_{\lambda,x}^*}\!\!
\left[
\left\|
\mathbf u-\mathbf u^\star
\right\| \mathbf{1}_{\mathcal E_R}\right] + \mathbb E_{Q_{\lambda,x}^*}\!\!
\left[
\left\|
\mathbf u-\mathbf u^\star
\right\| \mathbf{1}_{\mathcal E_R^c}\right]
\nonumber\\
&\leq
R
+
2\sqrt{mN} u_{\max} \,
Q_{\lambda,x}^*\!\left(
\mathcal E_R^c
\right)
\nonumber\\
&\overset{\eqref{eq:gibbs_tail_numerator_bound}}{\leq}
R
+
\frac{2\sqrt{mN} u_{\max}}{
p_{N,c}^{\mathrm{adm}}(r)
}
\exp\!\left(
-\frac{g_{N,c}(r,R)}{\lambda}
\right) \nonumber
\end{align}
which shows \eqref{eq:ideal_mppi_distance_to_global_optimum}.

\textbf{(3)} The upper bound \eqref{eq:ideal_mppi_distance_to_global_optimum} can be translated into \eqref{eq:nonconvex_ideal_mppi_suboptimality} using Assumption \ref{ass:local_lipschitz_cost}.

\textbf{(4)} 
Finally, the auxiliary radii may be chosen to shrink
with $\lambda$. Due to Assumption \ref{ass:unique_isolated_minimizer}, we can take any sequence $R_k\downarrow0$ and choose $r_k\in(0,R_k)$ such that $g_k:=g_{N,c}(r_k,R_k)>0$.
Since $p_{N,c}^{\mathrm{adm}}(r_k)>0$, we can select
$\lambda_k\downarrow0$ sufficiently fast that
\[
\frac{2\sqrt{mN}u_{\max}}
{p_{N,c}^{\mathrm{adm}}(r_k)}
e^{-g_k/\lambda_k}\le R_k.
\]
For $\lambda\in(\lambda_{k+1},\lambda_k]$,
using $(r_k,R_k)$ in the defining infimum gives
$\varepsilon_{N,c}(\lambda)\le 2R_k$.
Thus, although the separation gap and prior coverage
may decrease as the radii shrink, the temperature
can decrease sufficiently fast to compensate,
yielding $\varepsilon_{N,c}(\lambda) \rightarrow 0$ as $\lambda \rightarrow 0$.

\section*{Proof of Proposition \ref{prop:mppi_control_mismatch_hp}}
For any $x$ and prior distribution $P$, define for any $\tau>0$
\begin{equation}
    Q_{\lambda,x,\tau}(d \mathbf u)
    =
    \frac{
        \exp(-\tau J_N(x, \mathbf u)/\lambda)
    }{
        \mathbb E_{P}\!
    \left[
        \exp(-\tau J_N(x, \mathbf u)/\lambda)
    \right]
    }
    P(d \mathbf u).
\end{equation}
We first show that \begin{align}
    & \mathbb E_{P}\!\!
    \left[
        \left\|
            \mathbf u_{M,\lambda}(x)
            \!-\!
            \mathbf u_{\infty,\lambda}(x)
        \right\|^2
    \right] \!\leq\! \frac{16 mN u_{\max}^2 }{M} [\chi^2(Q_{\lambda,x}^* \| P)\!+\!1]  =
    \frac{16 mN u_{\max}^2}{M} \,
        \exp(\Gamma_\lambda(x)).
    \label{eq:gibbs_variance_bound}
\end{align}
And if in addition Assumption \ref{ass:convex_bounded_cost} holds, then 
   \begin{equation}
    \begin{split}
    \mathbb E_{P} &
    \left[
        \left\|
            \mathbf u_{M,\lambda}(x)
            -
            \mathbf u_{\infty,\lambda}(x)
        \right\|^2
    \right]  \leq
    \frac{16 mN u_{\max}^2}{M } 
        \exp\!\left(\frac{C_N(x)^2}{\lambda^2} \right).
    \end{split}
    \label{eq:gibbs_variance_bound_boundedJ}
\end{equation}
Define
\[
\begin{aligned}
Z_{\lambda,x}
&:=
\mathbb E_P\!\left[
\exp\!\left(-\frac{J_N(x,\mathbf u)}{\lambda}\right)
\right], \\
    \overline r_M
    &:=
    \frac{1}{M}
    \sum\nolimits_{m=1}^M
    \exp\!\left(-J_N(x,\mathbf{u}^{(m)})/\lambda\right) \big/ Z_{\lambda,x}, \\
    \overline Y_M
    &:=
    \frac{1}{M}
    \sum\nolimits_{m=1}^M
    \frac{\exp\!\left(-J_N(x,\mathbf{u}^{(m)})/\lambda\right)}{Z_{\lambda,x}}
    \left(
        \mathbf{u}^{(m)}-\mathbf u_{\infty,\lambda}
    \right).
    \end{aligned}
\]
Then
\begin{equation}
    \mathbf u_{M,\lambda} - \mathbf u_{\infty,\lambda}
    =
    \frac{
        \overline Y_M
    }{
        \overline r_M
    }.
    \label{eq:proof_error_ratio}
\end{equation}

\textbf{(1)} The expectation of $\overline r_M$ under the prior distribution $P$ can be obtained below. Since $\mathbf{u}^{(m)}$ are i.i.d, so 
\begin{equation}
\begin{split}
    \mathbb E_P \left[\overline r_M  \right] = \mathbb E_P \left[\exp\!\left(-J_N(x,\mathbf{u}^{(m)})/\lambda\right) \big/ Z_{\lambda,x}  \right] =1.
\end{split} \label{eq:mean_rM}
\end{equation}
\begin{equation}
\begin{split}
    \operatorname{Var}\left[\overline r_M  \right] = \frac{1}{M}\left(\mathbb E_P \!\left[ \exp\!\left(\frac{-2J_N(x,\mathbf{u}^{(m)})}{\lambda}\right)  \right] \Big/ Z_{\lambda,x}^2 -1 \right).
\end{split} \label{eq:var_rM}
\end{equation}

\textbf{(2)} The expectation of $\overline Y_M$ under the prior distribution $P$ can be obtained as below. Since
\(
\mathbf u^{(1)},\ldots,\mathbf u^{(M)}
\overset{\mathrm{i.i.d.}}{\sim} P,
\)
we have
\[
\begin{aligned}
&\mathbb E_P[\overline Y_M]
=
\mathbb E_P\!\left[
\exp\!\left(-\tfrac{J_N(x,\mathbf u)}{\lambda}\right)
\left(
\mathbf u-\mathbf u_{\infty,\lambda}
\right)
\right] \big/ Z_{\lambda,x}
\\
&=
\frac{\mathbb E_P\!\left[
\exp\!\left(-\tfrac{J_N(x,\mathbf u)}{\lambda}\right)
\mathbf u
\right]}{Z_{\lambda,x}}
-
\frac{\mathbf u_{\infty,\lambda}
\mathbb E_P\!\left[
\exp\!\left(-\tfrac{J_N(x,\mathbf u)}{\lambda}\right)
\right]}{Z_{\lambda,x}}.
\end{aligned}
\]
Since
\[
Q_{\lambda,x}^*(d\mathbf u)
=
\frac{
\exp(-J_N(x,\mathbf u)/\lambda)
}{
Z_{\lambda,x}
}
P(d\mathbf u),
\]
we obtain
\[
\mathbb E_P\!\left[
\exp\!\left(-\frac{J_N(x,\mathbf u)}{\lambda}\right)
\mathbf u
\right]
=
Z_{\lambda,x}\,\mathbb E_{Q_{\lambda,x}^*}[\mathbf u]
=
Z_{\lambda,x}\,\mathbf u_{\infty,\lambda}.
\]
Substituting back into the equation of $\mathbb E[\overline Y_M]$ above yields
\begin{equation}
\mathbb E[\overline Y_M]=0.
\label{eq:E_Y_M}
\end{equation}

\textbf{(3)} Consider the event
\begin{equation}
    \mathcal A
    :=
    \left\{
        \overline r_M\geq\frac12
    \right\},
\end{equation}
On which
\(
    \frac{1}{\overline r_M^2}
    \leq 4.
\)
Hence
\[
\begin{aligned}
    &\mathbb E_P
    \left[
        \left\|
            \mathbf u_{M,\lambda} - \mathbf u_{\infty,\lambda}
        \right\|^2
        \mathbf 1_{\mathcal A}
    \right] \\
    &\leq
    4 \,
    \mathbb E_P
    \left[
        \|\overline Y_M\|^2
    \right]\\
    & = \frac{4}{Z_{\lambda,x}^2}\mathbb E_P\!\! \left[ \left\| \frac{1}{M} \!
    \sum_{m=1}^M \!
    \exp\!\!\left(\tfrac{-J_N(x,\mathbf{u}^{(m)})}{\lambda}\right)\!
    \left(
        \mathbf{u}^{(m)}\!-\!\mathbf u_{\infty,\lambda}
    \right) \right\|^2_2 \right] \\
    &=
    \frac{4}{MZ_{\lambda,x}^2}
\mathbb E_P
\left[
    \exp\!\left(
        -\frac{2J_N(x,\mathbf u)}{\lambda}
    \right)
    \left\|
        \mathbf u-\mathbf u_{\infty,\lambda}
    \right\|^2
\right]
    \end{aligned}
\]
where the last equation holds because $\mathbf u^{(m)}$ are independent for all $m=1,2,...,M$ and \eqref{eq:E_Y_M}. In addition, since \(
\|\mathbf u^{(m)}\|_\infty\leq u_{\max},
\)
and
\(
\mathbf u_{\infty,\lambda}
=
\mathbb E_{Q_{\lambda,x}^*}[\mathbf u],
\) so
\begin{equation}
\begin{aligned}
    &\mathbb E_P
    \left[
        \left\|
            \mathbf u_{M,\lambda} - \mathbf u_{\infty,\lambda}
        \right\|^2
        \mathbf 1_{\mathcal A}
    \right] \\
    &\leq
    \frac{4}{MZ_{\lambda,x}^2}
\mathbb E_P
\left[
    \exp\!\left(
        -\frac{2J_N(x,\mathbf u)}{\lambda}
    \right)
    \left\|
        \mathbf u-\mathbf u_{\infty,\lambda}
    \right\|^2
\right] \\
& \leq \frac{16 mN u_{\max}^2}{M}
\mathbb E_P
\left[
    \exp\!\left(
        -\frac{2J_N(x,\mathbf u)}{\lambda}
    \right)
\right] \Big/ Z_{\lambda,x}^2 
    \end{aligned} \label{eq:u_variance_bound_eventA}
\end{equation}

\textbf{(4)} Now we control the complement event 
\begin{equation}
    \mathcal A^c
    :=
    \left\{
        \overline r_M <\frac12
    \right\}.
\end{equation}
We have 
\begin{equation}
\begin{split}
    \mathbb P(\mathcal A^c) &= \mathbb P \left(\overline r_M <\frac12 \right) \leq \mathbb P \left(\left|\overline r_M -1 \right| \geq\frac12 \right) \\
    & = \mathbb P \left(\left|\overline r_M - \mathbb E_P[\overline r_M] \right| \geq\frac12 \right) \\
    & \leq 4 \operatorname{Var}(\overline r_M) \\
    & \overset{\eqref{eq:var_rM}}{=}\! \frac{4}{M} \!\left( \mathbb E_P \!\left[ \exp\!\left(\frac{-2J_N(x,\mathbf{u})}{\lambda}\right)  \right] \Big/ Z_{\lambda,x}^2 -1 \!\right)
\end{split} \label{eq:Ac_event_probability}
\end{equation}
where the second inequality comes from Chebyshev's inequality. Then 
\begin{equation}
\begin{aligned}
    &\mathbb E_P
    \left[
        \left\|
            \mathbf u_{M,\lambda} - \mathbf u_{\infty,\lambda}
        \right\|^2
        \mathbf 1_{\mathcal A^c}
    \right] \\
    &\leq
    4 mN u_{\max}^2 \mathbb P(\mathcal A^c) \\
    & \overset{\eqref{eq:Ac_event_probability}}{\leq} \frac{16 mN u_{\max}^2}{M} \,\mathbb E_P \!\left[ \exp\!\left(\frac{-2J_N(x,\mathbf{u})}{\lambda}\right)  \right] \Big/ Z_{\lambda,x}^2
    \end{aligned} \label{eq:u_variance_bound_eventAc}
\end{equation}
where the first inequality comes from bounded input assumption.

\textbf{(5)} Finally, decomposing the expectation over $\mathcal A$ and $\mathcal A^c$ and applying
\eqref{eq:u_variance_bound_eventA} and
\eqref{eq:u_variance_bound_eventAc},
\begin{align}
    &
    \mathbb E
    \left[
        \left\|
            \mathbf u_{M,\lambda} - \mathbf u_{\infty,\lambda}
        \right\|^2
    \right]
    \notag\\
    &=
    \mathbb E
    \left[
        \left\|
            \mathbf u_{M,\lambda} - \mathbf u_{\infty,\lambda}
        \right\|^2
        \mathbf 1_{\mathcal A}
    \right]
    +
    \mathbb E
    \left[
        \left\|
            \mathbf u_{M,\lambda} - \mathbf u_{\infty,\lambda}
        \right\|^2
        \mathbf 1_{\mathcal A^c}
    \right]
    \notag \\
    &\leq
    \frac{16 mN u_{\max}^2}{M} \,\mathbb E_P \!\left[ \exp\!\left(\frac{-2J_N(x,\mathbf{u})}{\lambda}\right)  \right] \Big/ Z_{\lambda,x}^2 \notag\\
    & \overset{\eqref{eq:var_rM_chi_divergence}}{\leq } \frac{16 mN u_{\max}^2 }{M} [\chi^2(Q_{\lambda,x}^* \| P) +1] \notag\\
    & \overset{\eqref{eq:var_rM_var_J}}{=}\frac{16 mN u_{\max}^2}{M}
        \,\exp(\Gamma_\lambda(x))
\end{align}
where the last two formulas come from Lemma \ref{lem:var_rM_var_J}, which proves \eqref{eq:gibbs_variance_bound}. To show \eqref{eq:gibbs_variance_bound_boundedJ}, we use the upper bound on $J(\cdot,\cdot)$ from Assumption \ref{ass:convex_bounded_cost} and Popoviciu's inequality to obtain
\[
\operatorname{Var}_{Q}
    \left(
        J_N(x,\mathbf u)
    \right) \leq C_N(x)^2
\]
for any distribution $Q$. Substituting this into \eqref{eq:Gamma_t} and \eqref{eq:gibbs_variance_bound} yields the result.

We now proceed to show \eqref{eq:mppi_control_mismatch_hp}. Applying Markov's inequality to \eqref{eq:gibbs_variance_bound},\eqref{eq:gibbs_variance_bound_boundedJ} and setting $$\frac{16 mN u_{\max}^2}{M \epsilon^2} \,
        \exp(\Gamma_\lambda(x))=\delta, \ \frac{16mNu_{\max}^2}
{M\varepsilon^2} \,
    \exp\left(
        \frac{C_N(x)^2}{\lambda^2}
    \right)=\delta$$ respectively yields the conclusion.

\section*{Proof of Proposition
\ref{prop:mppi_finite_sample_drift}}
Fix any $x\in\Omega_c$ and $\delta\in(0,1)$ satisfying the conditions of the proposition. By Proposition~
\ref{prop:mppi_control_mismatch_hp},
the event
\[
\mathcal E_{x,\delta}
:=
\left\{
\|\mathbf u_{M,\lambda}(x)
-\mathbf u_{\infty,\lambda}(x)\|
\leq\eta_{M,\lambda,N}(x,\delta)
\right\}
\]
satisfies
\(
\mathbb P(\mathcal E_{x,\delta}\mid x)
\geq1-\delta.
\)
On $\mathcal E_{x,\delta}$, by the premise of Proposition~\ref{prop:suboptimality_ideal_mppi} and Assumption~\ref{ass:local_lipschitz_cost},
\begin{align}
J_N(x,\mathbf u_{M,\lambda}(x))
-V_N(x) &\leq
L_{c,N}^u
\|\mathbf u_{M,\lambda}(x)
-\mathbf u^\star(x)\|
\nonumber\\
&\leq L_{c,N}^u\Big[
\varepsilon_{N,c}(\lambda)
+\eta_{M,\lambda,N}(x,\delta)
\Big]
\nonumber\\
&=\Delta_{N,c}^{M}(x,\lambda,\delta)
\leq\tau.
\label{eq:hp_drift_cost_gap}
\end{align}
For brevity, denote the realized
suboptimality gap by
\(
e:=
J_N(x,\mathbf u_{M,\lambda}(x))
-V_N(x)
\).

Let
$\{\hat x_k\}_{k=0}^{N}$
be the predicted trajectory generated
by $\mathbf u_{M,\lambda}(x)$, and write
\(
\ell_k:=
\ell(\hat x_k,u_{M,k}),
\, k=0,\ldots,N-1.
\)
Since $\ell(x,u)\geq\sigma(x)$
and the stage costs are nonnegative, we have
\(
\sum_{k=0}^{N-1}\sigma(\hat x_k)
\leq
J_N(x,\mathbf u_{M,\lambda}(x))
=V_N(x)+e
\leq c+\tau.\)
Consequently,
\(
\sigma(\hat x_k)\leq c+\tau,
\, \forall k=0,\ldots,N-1.
\)
By Assumption~\ref{ass:no_terminal_cost_condition} and the definition of $\gamma_{c,\tau}$, we have
\begin{equation}
V_j(\hat x_k)
\leq\gamma_{c,\tau}\sigma(\hat x_k),
\,
j\geq1,\, \forall k=0,\ldots,N-1.
\label{eq:hp_drift_local_controllability}
\end{equation}
Applying the suboptimal MPC decrease argument of \citep[Theorem 1]{fiedler2025analysis} to the region $\{x\in \mathbb R^n \mid \sigma(x) \leq c+\tau\}$ with the constant $\gamma_{c,\tau}$ yields
\begin{align}
&V_N(f(x,\kappa_{M,\lambda}(x)))
-V_N(x) \leq
-\alpha_{N,c,\tau}\sigma(x)
+\gamma_{c,\tau}e \leq
-\alpha_{N,c,\tau}\sigma(x)
+\gamma_{c,\tau}
\Delta_{N,c}^{M}(x,\lambda,\delta),
\end{align}
where the last inequality follows from \eqref{eq:hp_drift_cost_gap}. This inequality holds on $\mathcal E_{x,\delta}$, whose conditional probability is at least $1-\delta$. This completes the proof.

\section*{Proof of Proposition \ref{thm:mppi_semi-global_recurrence}}

The proof proceeds in three steps. First, we select the prediction horizon $N$, temperature $\lambda$, and sample size $M$ to ensure that the high-probability Lyapunov decrease condition in Proposition~\ref{prop:mppi_finite_sample_drift} holds with a sufficiently small residual. Second, we employ a bootstrap argument to establish finite-time invariance of some compact set with high probability. Third, we show practical recurrence with high probability. Specifically, the Lyapunov decrease condition ensures that, whenever the state lies in a prescribed compact sublevel set, it remains in that set at the next time step on a high-probability event. Repeating this argument yields finite-time invariance with probability. On the same event, the Lyapunov function decreases strictly outside the practical recurrence set,
which guarantees entry into that set within finite time.

\textit{(1) Parameter selection.}

By Assumption~\ref{ass:no_terminal_cost_condition},
\(
\sigma(x)\leq V_N(x)\leq\bar a(x)\sigma(x),
\, \forall N\in\mathbb N.
\)
Since $\bar a$ is locally bounded, we can choose $c>0$ such that
\(
\sup_{\|x\|\leq r}\bar a(x)\sigma(x)<c.
\)
Consequently,
\(
V_N(x_0)\leq c,
\, \forall\|x_0\|\leq r.
\)

Fix any $\tau>0$ and define
\(
\gamma_{c,\tau}
:=
\sup_{\sigma(x)\leq c+\tau}\bar a(x)<\infty.
\)
Since $\gamma_{c,\tau}$ is independent of $N$, there exists $\underline N\geq2$
such that
\[
\alpha_{N,c,\tau}
=
1-(\gamma_{c,\tau}-1)^2
\left(
\frac{\gamma_{c,\tau}-1}{\gamma_{c,\tau}}
\right)^{N-2}>0
\]
for all $N\geq\underline N$. Fix any such $N\geq\underline N$ and define
\begin{equation}
T:=\left\lfloor\frac{c}{\varepsilon}
\right\rfloor+1,
\,
\delta:=\frac{p}{2T}.
\label{eq:recurrence_T}
\end{equation}
For brevity, write
$
\alpha:=\alpha_{N,c,\tau},
\,
\gamma:=\gamma_{c,\tau},
$
and define
\begin{equation}
\theta:=
\min\left\{
\tau,\,
\frac{\varepsilon}{\gamma},\,
\frac{\alpha c}{\gamma^2}
\right\}>0.
\label{eq:recurrence_design_tolerance}
\end{equation}
By Proposition~
\ref{prop:suboptimality_ideal_mppi},
$\varepsilon_{N,c}(\lambda)\to0$
as $\lambda\to0$. Hence, there exists
$\bar\lambda_N>0$ such that
\[
L_{c,N}^u\varepsilon_{N,c}(\lambda)
<\frac{\theta}{2},
\,
\forall\lambda\in(0,\bar\lambda_N).
\]
For any such $\lambda$, define
\[
C_{N,c}:=
\max_{\substack{x\in\Omega_c\\
\mathbf u\in\mathbb U^N}}
J_N(x,\mathbf u)<\infty.
\]
By Proposition~
\ref{prop:mppi_control_mismatch_hp},
there exists
$\underline M(N,\lambda,\delta)\in\mathbb N$
such that, for all $M\geq\underline M$,
\[
4L_{c,N}^u u_{\max}
\sqrt{\frac{mN}{M\delta}}
\exp\left(
\frac{C_{N,c}^2}{2\lambda^2}
\right)
<\frac{\theta}{2}.
\]
Therefore, by the definition of
$\Delta_{N,c}^{M}$ in
Proposition~\ref{prop:mppi_finite_sample_drift},
\begin{equation}
\sup_{x\in\Omega_c}
\Delta_{N,c}^{M}(x,\lambda,\delta)
\leq\theta.
\label{eq:recurrence_uniform_hp_residual}
\end{equation}
In particular, due to our design in \eqref{eq:recurrence_design_tolerance}, we have
$\Delta_{N,c}^{M}(x,\lambda,\delta) \leq\tau$ for all $x\in\Omega_c$, so Proposition~
\ref{prop:mppi_finite_sample_drift} is applicable throughout $\Omega_c$.

\textit{(2) High-probability invariance.}

By Proposition~
\ref{prop:mppi_finite_sample_drift}
and \eqref{eq:recurrence_uniform_hp_residual},
for every $x_t\in\Omega_c$,
\begin{equation}
\begin{aligned}
\mathbb P\Big(
V_N(x_{t+1})-V_N(x_t)
\leq-\alpha\sigma(x_t)
+\gamma\theta
\mid x_t
\Big)\geq1-\delta.
\end{aligned}
\label{eq:recurrence_hp_one_step}
\end{equation}

We first utilize \eqref{eq:recurrence_hp_one_step} to establish finite-time containment of the closed-loop trajectory in $\Omega_c$ with high probability. Since $V_N(x)\leq\gamma\sigma(x)$ for all $x\in\Omega_c$ according to Assumption \ref{ass:no_terminal_cost_condition}, the decrease inequality in \eqref{eq:recurrence_hp_one_step} implies
\begin{align}
V_N(x_{t+1})
&\leq V_N(x_t)-\alpha\sigma(x_t)
+\gamma\theta
 \leq
\left(1-\frac{\alpha}{\gamma}\right)
V_N(x_t)+\gamma\theta.
\label{eq:recurrence_hp_contraction}
\end{align}
By the choice of $\theta$ in
\eqref{eq:recurrence_design_tolerance},
\(
\gamma\theta\leq\frac{\alpha c}{\gamma}.
\)
So whenever $x_t\in\Omega_c$ and \eqref{eq:recurrence_hp_contraction} holds,
\begin{align}
V_N(x_{t+1})
&\leq
\left(1-\frac{\alpha}{\gamma}\right)c
+\frac{\alpha c}{\gamma}
 =c.
\end{align}
Thus, on each successful sampling event,
\begin{equation}
x_t\in\Omega_c
\quad\Longrightarrow\quad
x_{t+1}\in\Omega_c.
\label{eq:recurrence_one_step_containment}
\end{equation}

Let $\mathcal E_t$ denote the event that
the Lyapunov decrease inequality in
\eqref{eq:recurrence_hp_one_step}
holds at time $t$, and define
\[
\mathcal E:=
\bigcap_{t=0}^{T-1}\mathcal E_t.
\]
Since $x_0\in\Omega_c$, \eqref{eq:recurrence_one_step_containment} implies, by induction, that
$
x_t\in\Omega_c,\, t=0,\ldots,T,
$
on $\mathcal E$. Since Algorithm~\ref{alg:mppi} draws fresh independent samples at each step, Proposition~\ref{prop:mppi_finite_sample_drift} ensures that
\(
\mathbb P(\mathcal E_t\mid x_t)
\geq1-\delta
\). Therefore, a union bound argument over the T steps gives us
\begin{align}
\mathbb P(\mathcal E)
&\geq(1-\delta)^T
 \geq1-T\delta
 \overset{\eqref{eq:recurrence_T}}{=}1-\frac p2>1-p.
\label{eq:recurrence_joint_probability}
\end{align}

\textit{(3) High-probability practical recurrence.}
We now establish practical recurrence on the same event $\mathcal E$. Recall the practical recurrence set
$
\mathcal A_{c,\varepsilon}
:=
\left\{
x\in\Omega_c:
\alpha_{N,c,\tau}\sigma(x)
\leq2\varepsilon
\right\}.
$
Whenever $x_t\notin \mathcal A_{c,\varepsilon}$, the established containment implies
$x_t\in\Omega_c$, and hence
\(
\alpha\sigma(x_t)>2\varepsilon.
\)
Since $\gamma\theta\leq\varepsilon$ due to \eqref{eq:recurrence_design_tolerance},
the decrease inequality in
\eqref{eq:recurrence_hp_one_step}
gives
\begin{equation}
V_N(x_{t+1})-V_N(x_t)
<-\varepsilon.
\label{eq:recurrence_strict_decrease}
\end{equation}
Suppose, for contradiction, that the trajectory does not enter $\mathcal A_{c,\varepsilon}$ during $t=0,\ldots,T$ on $\mathcal E$. Then \eqref{eq:recurrence_strict_decrease} holds for every $t=0,\ldots,T-1$. Summing these inequalities yields
\begin{align}
V_N(x_T)<V_N(x_0)-T\varepsilon
 \leq c-T\varepsilon
 <0,
\end{align}
where the last inequality follows from the choice of $T$ in \eqref{eq:recurrence_T}. This contradicts $V_N(x_T)\geq0$. Therefore, on $\mathcal E$, there exists $k\in\{0,\ldots,T\}$ such that $x_k\in\mathcal A_{c,\varepsilon}$. Combining this result with \eqref{eq:recurrence_joint_probability} gives
$
\mathbb P\left[
\exists k\in\{0,\ldots,T\}:
x_k\in\mathcal A_{c,\varepsilon}
\right]>1-p
$
for every $\|x_0\|\leq r$. This completes the proof.

\end{document}